\documentclass[11pt]{article}
\usepackage{amssymb,amsfonts,amsthm}

\newcommand{\be}[1]{\begin{equation}\label{#1}}
\newcommand{\ee}{\end{equation}}
\newcommand{\topp}[1]{\left\lceil{#1}\right\rceil}
\newcommand{\bott}[1]{\lfloor{#1}\rfloor}
\newcommand{\aaa}{{\cal A}}
\newcommand{\dd}{{\cal D}}
\newcommand{\kk}{{\cal K}}
\newcommand{\mm}{{\cal M}}
\newcommand{\nn}{{\cal N}}
\newcommand{\pp}{{\cal P}}
\newcommand{\sss}{{\cal S}}
\newcommand{\gr}{{\cal G}}
\newcommand{\fff}{{\cal F}}
\newcommand{\barkk}{\overline{\kk}}
\newcommand{\ve}{\varepsilon}
\newcommand{\la}{\langle\,}
\newcommand{\ra}{\,\rangle}
\newcommand{\Ker}{\mathop{\rm Ker}}
\newcommand{\Sq}{\mathop{\rm Sq}}

\newcommand{\E}{{\bf E}}
\newcommand{\F}{{\bf F}}
\newcommand{\hh}{{\cal H}}
\newcommand{\wtk}{\widetilde\kk}
\newcommand{\wtD}{\widetilde\Delta}

\newcommand{\te}{{\tilde e}}
\newcommand{\tf}{{\tilde f}}
\newcommand{\tp}{{\tilde p}}
\newcommand{\tq}{{\tilde q}}
\newcommand{\tr}{{\tilde r}}

\newcommand{\tu}{{\tilde u}}
\newcommand{\tv}{{\tilde v}}
\newcommand{\tx}{{\tilde x}}

\newcommand{\bfc}{{\bf H}}
\newcommand{\plf}{{\hbox{\bf PLF}}}
\newcommand{\PLF}{\mathop{\rm PLF}}
\newcommand{\Homeo}{\mathop{\rm Homeo}}

\begin{document}

\newcounter{ppp}
\setcounter{ppp}{1}

\title
{Diagram groups are totally orderable}
\author{V.\,S.\,Guba,\ \,M.\,V.\,Sapir\thanks{The research of the first
author was partially supported by the RFFI grant 99--01--00894 and the
INTAS grant 99--1224. The research of the second author was supported in
part by the NSF grants DMS 0072307 and 9978802 and the US-Israeli BSF grant
1999298.}}
\date{}

\maketitle

\begin{abstract}

In this paper, we introduce the concept of the independence graph
of a directed $2$-complex. We show that the class of diagram
groups is closed under graph products over independence graphs of
rooted $2$-trees. This allows us to show that a diagram group
containing all countable diagram groups is a semi-direct product
of a partially commutative group and R. \,Thompson's group $F$. As
a result, we prove that all diagram groups are totally orderable.
\end{abstract}

\theoremstyle{plain}
\newtheorem{thm}{Theorem}[section]
\newtheorem{lm}[thm]{Lemma}
\newtheorem{cy}[thm]{Corollary}
\newtheorem{df}[thm]{Definition}
\newtheorem{rk}[thm]{Remark}
\newtheorem{ex}[thm]{Example}
\newtheorem{prop}[thm]{Proposition}
\newtheorem{prob}[thm]{Problem}

\tableofcontents

\section{Introduction}

This paper is based on \cite{GuSa02a}, where we define many of the
concepts used below.

We say that $G$ is {\em left-orderable\/} ({\em right-orderable\/})
whenever there exists a total order $<$ on $G$ invariant under left
(right) multiplication, that is, $a<b$ implies $ca<cb$ ($ac<bc$) for any
$a,b,c\in G$. It is easy to see that a group is left-orderable if and
only if it is right-orderable. The group $G$ is {\em orderable\/}
whenever there exists a total order $<$ on $G$ such that $a<b$ implies
$ca<cb$ and $ac<bc$ for any $a,b,c\in G$. One of the main results of this
paper is that all diagram groups are totally orderable. This answers
in positive \cite[Problem 17.6]{GuSa97}. When this paper was almost
completed, we found out from B.\,Wiest that he had recently proved
that diagram groups are left-orderable (which is weaker than being
orderable). First he proved that diagram groups are embeddable into
a certain braid group on infinitely many strings, then he proved that
this braid group is left-orderable \cite{Wi02}. Note that this very
interesting approach is completely independent from ours.
\vspace{1ex}

We use standard notation for conjugation and commutators. If $a$,
$b$ are elements of a group $G$, then $a^b=b^{-1}ab$,
$[a,b]=a^{-1}a^b=a^{-1}b^{-1}ab$. If $A$, $B$ are subgroups of
$G$, then $[A,B]$ is the subgroup of $G$ generated by all
commutators of the form $[a,b]$, where $a\in A$, $b\in B$.

A group $G$ is called {\em partially commutative\/} ({\em right angled
Artin group\/} or {\em graph group\/}) if it has a group presentation
such that all defining relations have the form $[a,b]=1$, where $a$, $b$
are generators.

We will use some properties of the R.\,Thompson group $F$ which is the
diagram group of the Dunce hat $\hh_0=\la x\mid x^2=x\ra$ (see
\cite{GuSa02,GuSa02a}). Recall \cite{CFP} that $F$ is isomorphic to the
group of all increasing continuous piecewise linear functions from the
interval $[0,1]$ onto itself such that all singularities (breakpoints of
the derivative) occur only at finitely many dyadic rational points (points
of the form $m/2^n$) and all slopes are integer powers of $2$. The group
operation is the composition of functions (we shall write function symbols
to the right of the argument).

In \cite{GuSa02a}, we proved that the diagram group $\gr_1=\dd(\hh_1,x)$
of the directed $2$-complex $\hh_1=\la x\mid x^2=x,x=x\ra$ is
{\em universal\/}, that is, it contains all countable diagram groups. Notice
that $\hh_1$ is obtained from the Dunce hat by adding one $2$-cell of the
form $x=x$.
\vspace{0.5ex}

Throughout the paper, we shall always assume that when we add a $2$-cell,
we add its inverse as well.
\vspace{1ex}

In this paper, we first prove, using diagram products of groups
\cite{GuSa99}, that if we add $2$-cells of the form $e=e$ to a
directed $2$-complex $\kk$, then the diagram groups of the new
complex are semi-direct products of partially commutative groups
and the diagram groups of the initial complex. Moreover, the
partially commutative group can be explicitly described as a
diagram group of some directed $2$-complex. In order to do that,
we use the concept of a rooted $2$-tree from \cite{GuSa02a} and
introduce the concept of the independence graph of a directed
$2$-complex. We show that the class of diagram groups is closed
under graph products over independence graphs of rooted $2$-trees
(in particular, it is closed under countable direct products). We
also use natural representations of diagram groups and groupoids
by homeomorphisms of intervals of the real line using transition
functions as in \cite{GuSa97}. As a result, we obtain a complete
description of the structure of $\gr_1$, and prove that $\gr_1$
is orderable. That, in turn, implies that all diagram groups are
totally orderable because the orderability is a local property.

\section{Diagram products of groups and graph products of diagram groups}
\label{dps}

Let us recall the definition of the diagram product of groups (we
translate the definition from \cite{GuSa99} into the language of
directed $2$-complexes). Let $\kk$ be a directed $2$-complex with
the set of edges $\E$ and the set of positive $2$-cells $\F^+$. Let
$\gr_E=\{\,G_e\mid e\in\E\,\}$ be a collection of groups. Finally, let
$w$ be a nonempty $1$-path in $\kk$. Then the diagram product
$\dd(\gr_E;\kk,w)$ of $\gr_E$ over $\kk$ with base $w$ is the
fundamental group with base $w$ of the following $2$-complex $C(\gr_E)$
of groups \cite{Bri} described below:

\begin{itemize}

\item the underlying $2$-complex is the $2$-skeleton of the Squier
complex $\Sq(\kk)$,

\item the vertex group $G(s)$ assigned to a vertex $s=e_1\cdots e_n$
($e_i\in\E$) is the direct product $G_{e_1}\times\cdots\times G_{e_n}$,

\item the edge group assigned to the edge $x=(p,f,q)$ is the direct
product $G(p)\times G(q)$,

\item the embeddings $\iota_x$ of the edge group $G(p)\times G(q)$
into the vertex groups $G(p\topp{f}q)=G(p)\times G(\topp{f})\times
G(q)$ and $\tau_x$ of $G(p)\times G(q)$ into
$G(p\bott{f}q)=G(p)\times G(\bott{f})\times G(q)$ are coordinate-wise,

\item the group assigned to every $2$-cell in $\Sq(\kk)$ is trivial.
\end{itemize}

The procedure of calculating the diagram product can be described as
follows. First we choose a spanning forest ${\cal T}$ in $\Sq(\kk)$.
Then every edge $x$ of $\Sq(\kk,w)$ defines an element $\bar x$ of
the fundamental group $\pi_1(Sq(\kk),w)\cong\dd(\kk,w)$. Namely, this
is the element represented by the loop $p_{\iota(x)}^{-1}xp_{\tau(x)}$,
where $p_v$ is the geodesic path from $v$ to $w$ in ${\cal T}$.
According to formulas (4), (5) of \cite{GuSa99}, the diagram product
$\dd(\gr_E;\kk,w)$ is isomorphic to the group
\be{dpform}
\left.\mathop{\ast}\limits_{v}G(v)\ast\dd(\kk,w)\right/\nn,
\ee
where the free product of the groups $G(v)$ is taken over all
vertices of $\Sq(\kk,w)$ and $\nn$ is the normal closure of the
following defining relations:
\be{descr}
\iota_x(g)^{\bar x}=\tau_x(g)\mbox{  for every positive edge $x$ of }
\Sq(\kk,w), \mbox{ where }g\in G(x).
\ee

Theorem 4 of \cite{GuSa99} shows that this algebraic construction
corresponds to the following topological construction provided all
groups in $\gr_E$ are diagram groups.

The directed $2$-complex on Figure \theppp\ will be called a {\em
switch\/}. It has $4$ vertices, edges $a$, $b$, $c$, $s$ and one
positive $2$-cell $a=bsc$.

\begin{center}
\unitlength=0.5mm \special{em:linewidth 0.4pt}
\linethickness{0.4pt}
\begin{picture}(70.00,20.00)
\put(0.00,16.00){\vector(1,0){20.00}}
\put(20.00,16.00){\vector(1,0){30.00}}
\put(50.00,16.00){\vector(1,0){20.00}}
\put(37.00,6.00){\vector(1,0){00.00}}
\bezier{324}(0.00,16.00)(34.00,-4.00)(70.00,16.00)
\put(11.00,19.00){\makebox(0,0)[cc]{$b$}}
\put(35.00,19.00){\makebox(0,0)[cc]{$s$}}
\put(60.00,19.00){\makebox(0,0)[cc]{$c$}}
\put(35.00,1.00){\makebox(0,0)[cc]{$a$}}
\end{picture}

\nopagebreak[4] Figure \theppp.
\end{center}
\addtocounter{ppp}{1}

Let $\kk$ be a directed $2$-complex with the set of edges $\E$. For every
$e\in\E$, let $\kk_e$ be a directed $2$-complex with a distinguished
$1$-path $p_e$, where $G_e\cong\dd(\kk_e,p_e)$. For each edge $e\in\E$ we
attach a switch of the form $e=b_es_ec_e$ to $\kk$. Then attach $\kk_e$
by subdividing $s_e$ into $|p_e|$ edges and gluing $p_e$ with $s_e$.
As a result, we obtain a new directed $2$-complex $\barkk$. Note that
$G_e$ is also isomorphic to the diagram group $\dd(\hat\kk_e,e)$, where
$\hat\kk_e$ is $\kk_e$ with the switch added. This follows from
\cite{GuSa02a}, Theorem 4.1, part 1 and Corollary 3.4.

One can easily check that this description coincides with the
construction in \cite[Theorem 4]{GuSa99} if $\kk=\kk_\pp$ for some
semigroup presentation $\pp$. The only difference is that in
\cite{GuSa99} we used (in the new terminology) switches of the
form $a=bsb$ instead of $a=bsc$. This difference is insignificant
and the proof of \cite[Theorem 4]{GuSa99} can be easily
generalized to arbitrary directed $2$-complexes. We also need to
describe explicitly the isomorphism, which was presented in the
proof of that theorem. In order to do this, we need to introduce
some notation.

Let $v$ be a $1$-path in $\kk$ that contains an edge $z$. One can
write $v=v'zv''$ for some $1$-paths $v'$, $v''$. Then the vertex
group $G(v)$ in the complex $C(\gr_E)$ of groups equals the direct
product $G(v')\times G_z\times G(v'')$. The direct factor $G_z$ in
this direct product projects onto a subgroup in the diagram
product. This subgroup will be denoted by $G(v',z,v'')$. According
to (\ref{dpform}), all these subgroups together with the image of
diagram group $\dd(\kk,w)$ generate the diagram product.

For any vertex $v$ of $\Sq(\kk,w)$, we let $P_v$ be the
$(v,w)$-diagram over $\kk$ corresponding to the geodesic path
$p_v$ in ${\cal T}$ from $v$ to $w$. The next theorem immediately
follows from \cite[Theorem 4]{GuSa99} (more precisely, from the
straightforward generalization of this theorem to directed
$2$-complexes) and its proof.

\begin{thm}
\label{dpth} For every nonempty $1$-path $w$ in $\kk$, the diagram
product $\dd(\gr_E;\kk,w)$ is isomorphic to the diagram group
$\dd(\barkk,w)$. The isomorphism takes an element $g$ of the
subgroup $G(v',z,v'')$ to the diagram
$P_v^{-1}(\ve(v')+\Delta_g+\ve(v''))P_v$, where $v=v'zv''$,
$\Delta_g$ is a $(z,z)$-diagram that represents the element of
$\dd(\hat\kk_z,z)$ corresponding to $g$.
\end{thm}

This theorem gives an explicit representation by a diagram from the
group $\dd(\barkk,w)$ for any element of the diagram product
$\dd(\gr_E;\kk,w)$ . Formula (\ref{dpform}) says that the diagram
product is generated by the images of the diagram group $\dd(\kk,w)$
and the vertex groups $G(v)$. The elements of $\dd(\kk,w)$ are already
$(w,w)$-diagrams over $\barkk$. The vertex group $G(v)$ is a direct
product of the groups of the form $G(v',z,v'')$ so we know the
$(w,w)$-diagrams that represent elements of $G(v)$.

In this paper, we consider diagram products of groups in the case
when the underlying directed $2$-complex $\kk$ is a rooted
$2$-tree. Recall \cite[Theorem 7.2]{GuSa02a} that a rooted
$2$-tree with root $w$ is a directed $2$-complex $\kk$ satisfying
the following conditions:

\begin{description}
\item[T1] \ For any vertex $o$ of $\kk$, there exist a $1$-path
from $\iota(w)$ to $\tau(w)$ containing $o$.
\item[T2] \ Every two
$1$-paths in $\kk$ with the same endpoints are (directly)
homotopic in $\kk$.
\item[T3] \ The diagram group $\dd(\kk,w)$ is
trivial.
\end{description}

We are going to show now that if $\kk$ is a rooted $2$-tree then
any diagram product of groups over $\kk$ is a graph product of
these groups.

Let $\kk$ be a directed $2$-complex with a distinguished $1$-path $w$. Two
(different) edges $e_1$, $e_2$ of $\kk$ will be called {\em independent\/}
provided there exists a $1$-path in $\kk$ homotopic to $w$ that contains
both $e_1$ and $e_2$. Consider the following {\em independence graph\/}
$\Gamma(\kk,w)$. Its vertices are all edges of $\kk$ that belong to
$1$-paths homotopic to $w$. Two vertices are adjacent in the graph if and
only if they represent independent edges.

Let us recall the concept of a {\em graph product\/} of groups. Let $\Gamma$
be a graph. Consider a family of groups $\{\,G_v, v\in V\}$, where $V$ is
the set of vertices of $\Gamma$. The {\em graph product of groups
$\{\,G_v, v\in V\}$\/} over the graph $\Gamma$ is the factor-group of
the free product $*_{v\in V} G_v$ over the normal subgroup generated by
all commutators $[g,h]$ for all $g\in G_s$, $h\in G_t$, where $s$, $t$
are adjacent vertices in $\Gamma$.

Clearly, if $\Gamma$ is a complete graph, then the graph product coincides
with the direct product, and if $\Gamma$ is the graph with no edges, then
the graph product coincides with the free product. Partially commutative
groups are precisely the graph products of free (cyclic) groups.
\vspace{1ex}

\begin{lm}
\label{dptriv}
Let $\kk$ be a directed $2$-complex and let $w$ be its nonempty $1$-path.
Suppose that $\kk$ satisfies conditions {\rm T2} and {\rm T3}. Let $\E$ be
the set of edges of $\kk$, and $\gr_E=\{\,G_e\mid e\in\E\,\}$ be a family of
groups. Then the diagram product of the family of groups
$\{\,G_e\mid e\in\E\,\}$ over $\kk$ with base $w$ is isomorphic to the graph
product of the same family of groups over the independence graph
$\Gamma(\kk,w)$.
\end{lm}

\proof
Let $\sss=\Sq(\kk,w)$. The vertices of $\sss$ are all $1$-paths in
$\kk$ that are homotopic to $w$. Let $v=z_1\cdots z_m$ ($z_i\in\E$) be a
vertex of $\sss$. In the $2$-complex of groups $C(\gr_E)$, the vertex
group $G(v)$ is the direct product $G_{z_1}\times\cdots\times G_{z_m}$.
Since the diagram group $\dd(\kk,w)$ is trivial, the diagram product is
generated by the images of all the vertex groups $G(v)$, $v\in\sss$
according to (\ref{dpform}).

Let $z$ be any vertex of $\Gamma(\kk,w)$. This means that $z$ belongs
to some $1$-path $v$ homotopic to $w$, that is, $v=v'zv''$ for some
$1$-paths $v'$, $v''$. We are going to show that the image of the
group $G(v',z,v'')$ in the diagram product depends on $z$ only.

Suppose that $v_1'zv_1''$, $v_2'zv_2''$ are $1$-paths in $\kk$
homotopic to $w$. The paths $v_1'$, $v_2'$ have the same endpoints
so they are homotopic. Then there exists a $(v_1',v_2')$-diagram
$\Delta'$ over $\kk$. Analogously, there exists a $(v_1'',v_2'')$-diagram
$\Delta''$ over $\kk$. Each of these diagrams can be represented as a
concatenation of atomic diagrams. Thus in order to prove that the groups
$G(v_1',z,v_1'')$, $G(v_2',z,v_2'')$ are identified in the diagram product,
it suffices to consider the case when one of the diagrams $\Delta'$ or
$\Delta''$ is trivial and the other one is an atomic diagram. These cases
are symmetric, so without loss of generality let $\Delta'$ be trivial and
let $\Delta''$ be the atomic diagram corresponding to a positive edge
$x=(p,f,q)$ in $\Sq(\kk,w)$. The edge $z$ is contained in $p$, that is,
$p=p'zp''$. We need to compare the groups $G(p',z,p''\topp{f}q)$ and
$G(p',z,p''\bott{f}q)$. Notice that the edge $x$ represents the trivial
element of the fundamental group $\pi_1(\Sq(\kk),w)$ since the group is
trivial itself. So the element $\bar x$ in formula (\ref{descr}) is trivial.
The edge group $G(x)=G(p)\times G(q)$ embeds naturally into the vertex groups
$G(p\topp{f}q)=G(p)\times G(\topp{f}\times G(g)$ and
$G(p\bott{f}q)=G(p)\times G(\bott{f}\times G(g)$. These embeddings
have been denoted by $\iota_x$ and $\tau_x$. Formula (\ref{descr})
now says that $\iota_x(g)=\tau_x(g)$ for any $g\in G(x)$. Thus the
images of $G(p)$ in both vertex groups are identified in the diagram
product. In particular, the direct factors of the form $G_z$ are also
identified. But these are the groups $G(p',z,p''\topp{f}q)$ and
$G(p',z,p''\bott{f}q)$. So they coincide in the diagram product.

Thus the images of groups $G_z$, $z\in\E$, generate the diagram
product. The only defining relations we impose on the free product
of these groups, are commutativity relations involving elements of
different factors of the direct products $G(v)$, $v\in\Sq(\kk,w)$.
If $z_1$, $z_2$ are independent edges and $u_0z_1u_1z_2u_2$ is a
$1$-path homotopic to $w$, then the groups $G(z_1)$ and $G(z_2)$
pairwise commute in the diagram product because they are the factors
of the direct product
$G(u_0)\times G_{z_1}\times G(u_1)\times G_{z_2}\times G(u_2)$.
These are the only defining relations we have. So the diagram product
is the free product of the groups $G_z$, where $z$ occurs in a
$1$-path homotopic to $w$, factored by the relations of commutativity
of the form $[G_{z_1},G_{z_2}]=1$ for each pair of independent edges
$z_1$, $z_2$. This completes the proof.
\endproof

\begin{rk}
\label{expl}
{\rm
In case the groups $G_z$ are diagram groups, we are interested in
finding an explicit form of their elements in the diagram product.
Suppose that $g\in G_z$ is represented by a $(z,z)$-diagram $\Delta_z$
over $\hat\kk_z$. For any $1$-paths $v'$, $v''$ in $\kk$, where
$v'zv''$ is homotopic to $w$ in $\kk$, we can consider the diagram
of the form
\be{exfor}
P^{-1}(\ve(v')+\Delta_g+\ve(v''))P,
\ee
where $P$ is a $(v,w)$-diagram over $\kk$. If $P=P_v$ from the statement
of Theorem \ref{dpth}, then $g\in G(v',z,v'')$ is represented in
the diagram product by (\ref{exfor}). However, if the diagram group
$\dd(\kk,w)$ is trivial, then, replacing $P_v$ by $P$ leads to an
equivalent diagram. So in this case (\ref{exfor}) gives us the desired
representation of an element $g\in G_z$. Since the groups $G(v',z,v'')$
are independent of $v'$, $v''$, the diagram we get also depends on $z$
and $g\in G_z$ only.
}
\end{rk}

Let us call a graph $\Gamma$ {\em appropriate\/} if the class of
diagram groups is closed under graph products over $\Gamma$. If
every two $1$-paths in $\kk$ with the same endpoints are homotopic
and the diagram group $\dd(\kk,w)$ is trivial, then Lemma
\ref{dptriv} shows that the independence graph $\Gamma(\kk,w)$ is
appropriate. Clearly, the class of appropriate graphs is closed
under taking full subgraphs. (Indeed, we can assign a trivial
group to each vertex not in the subgraph.) On the other hand,
\cite[Theorem 30]{GuSa99} shows that an appropriate graph cannot
contain a cycle of odd length $\ge5$ as a full subgraph.

Lemmas \ref{dptriv} and \cite[Lemma 7.1]{GuSa02a} provide us with
a large class of appropriate graphs.

\begin{thm}
\label{approp}
Let $\kk$ be a rooted $2$-tree with the root $w$. Then the independence
graph $\Gamma(\kk,w)$ is appropriate.
\end{thm}

It is possible to prove that the converse of Theorem 2.4 also
holds: every appropriate graph is the independence graph of some
rooted $2$-tree. We shall include the proof in our future paper.

We know \cite{GuSa97,GuSa99} that the class of diagram groups is closed
under finite direct products and countable direct powers. Now we have a
stronger result.

\begin{thm}
\label{dirprod}
The class of diagram groups is closed under countable direct products.
\end{thm}

\proof
Direct products are graph products over complete graphs. Thus by Theorem
\ref{approp}, it suffices to construct a rooted $2$-tree whose independence
graph contains a countable complete subgraph. It is not difficult to
understand that this is true for any ``sufficiently branching" infinite
rooted $2$-tree.

In particular, let us start with the edge $e_0$ and then for each $n\ge0$
add a $2$-cell of the form $e_{2n}=e_{2n+1}e_{2n+2}$. In the resulting
rooted $2$-tree, the edges $e_1$, $e_3$, $e_5$, \dots with odd subscripts
are obviously pairwise independent.
\endproof

\section{Expansions of directed $2$-complexes}
\label{expan}

Let $\kk$ be a directed $2$-complex and let $\nu$ be a function from the
set of edges of $\kk$ to a set of (possibly infinite) cardinal
numbers. For any edge $e$ of $\kk$, let us add $\nu(e)$ positive
$2$-cells of the form $e=e$. The new directed $2$-complex $\kk_\nu$
is called an {\em expansion\/} of $\kk$. The new cells of the form
$e=e$ are called {\em leaves}. In particular, if $\kk=\hh_0$ is the
Dunce hat and $\nu$ is the function defined by $\nu(x)=1$, then
$\kk_\nu=\hh_1$ and $\dd(\kk_\nu,x)$ is the universal group $\gr_1$.

Expansions arise naturally when we consider directed $2$-complexes
$\kk$ with redundant $2$-cells. A $2$-cell $f$ is called {\em
redundant\/} if $\topp{f}$ is homotopic to $\bott{f}$ in the
complex $\kk\setminus\{f\}$ obtained from $\kk$ by removing
$f^{\pm 1}$. In this case, by \cite[Theorem 4.1, part 1]{GuSa02a},
we can add a new edge $e$ with $\iota(e)=\iota(f)$,
$\tau(e)=\tau(f)$ and a new $2$-cell of the form $\topp{f}=e$. The
diagram groups will not change. After that, by \cite[Theorem 4.1,
part 2]{GuSa02a}, we can replace the $2$-cell $f$ by the $2$-cell
$f'$ of the form $e=e$ without changing the diagram groups. Thus
adding redundant $2$-cells to a directed $2$-complex is
essentially equivalent to taking expansions.

Complexes with redundant $2$-cells are needed when we use complete
directed $2$-complexes in order to compute diagram groups. Recall that
we have quite powerful technical tools to compute diagram groups of
complete directed complexes \cite[Section 6]{GuSa02a}. If we want to
compute a diagram group of a directed $2$-complex, which is not necessarily
complete, then one way to do this is to embed $\kk$ into a bigger complex
$\kk'$, where $\kk'$ is complete. This can be done using a kind of the
Knuth -- Bendix completion procedure, so $\kk'$ is obtained from $\kk$ by
adding some redundant $2$-cells.

Notice that if $\kk'$ is obtained from $\kk$ by adding redundant
$2$-cells, then for every nonempty $1$-path $p$ of $\kk$, the
diagram group $H=\dd(\kk,p)$ is a retract of $G=\dd(\kk',p)$ (see
\cite[Lemma 4.1, part 1]{GuSa02a}) and the retraction can be described
explicitly. Thus if we know $H$, we can compute $G$.

The main goal of this section is to give a description of the diagram
groups of the expansion $\kk_\nu$ as semi-direct products of partially
commutative groups and the diagram groups of $\kk$.

Let $p$ be a nonempty $1$-path in $\kk$. Given a $(p,p)$-diagram $\Delta$
over $\kk_\nu$, we can collapse all cells of $\Delta$ that correspond to
$2$-cells of the form $e=e$ from $\kk_\nu\setminus\kk$. (Collapsing a cell,
means identifying its top and bottom path so the cell becomes an edge.) The
resulting diagram over $\kk$ is denoted by $\theta(\Delta)$. Clearly,
$\theta$ is a retraction from $\dd(\kk_\nu,p)$ to $\dd(\kk,p)$ because all
diagrams over $\kk$ are fixed by $\theta$.

Therefore, the group $\dd(\kk_\nu,p)$ is a semi-direct product of the kernel
$\aaa$ of $\phi$ and $G=\dd(\kk,p)$. We are going to prove that $\aaa$ is a
partially commutative group that can be presented as a diagram group of a
certain directed $2$-complex associated with $\kk_\nu$.
\vspace{0.5ex}

Let $p$ be any $1$-path in $\kk$. Recall \cite{GuSa02a} that a
{\em universal $2$-cover\/} of $\kk$ with base $p$ is a directed
$2$-complex $\mm$ over $\kk$ with a labelling map $\phi$ which
satisfy the following properties:

\begin{description}
\item[U1]\ $\mm$ is a rooted $2$-tree with root $\tp$, where
$\phi(\tp)=p$;
\item[U2]\ for any $1$-path $\tq$ in $\mm$ from $\iota(\tp)$ to
$\tau(\tp)$, the local map of $\mm$ at $\tq$ is bijective.
\end{description}

By \cite[Remark 8.2]{GuSa02a}, property U2 can be replaced by the
following property:

\begin{description}
\item[U2$'$]\ for any $1$-path $\tr$ in $\mm$ and for any $2$-cell
$f$ of $\kk$ such that $\topp{f}=\phi(\tr)$, there is exactly one
$2$-cell $\tf$ of $\mm$ labelled by $f$ with top path $\tr$.
\end{description}

Theorem \cite[Theorem 8.1]{GuSa02a} shows that for every directed
$2$-complex $\kk$ and a $1$-path $p$ in $\kk$, the universal
$2$-cover $\wtk_p$ exists and is unique up to an isomorphism that
preserves the labels.

Throughout this section, we fix a directed $2$-complex $\kk$, a $1$-path
$p$ in it, the expansion $\kk_\nu$, the retraction $\theta$, the kernel
$\aaa$ of $\theta$, the universal $2$-cover $\wtk_p$ with the root $\tp$,
and the labelling map $\phi$.

Let $\tilde\nu$ be the function on the edges of $\wtk_p$ that
sends every edge $\te$ to $\nu(\phi(\te))$. The corresponding
expansion of $\wtk_p$ is denoted by $\barkk_{\nu,p}$. There
exists a natural labelling morphism $\bar\phi$ from
$\barkk_{\nu,p}$ to $\kk_\nu$ that acts on $\wtk_p$ as $\phi$ and
acts as a one-to-one correspondence between the set of leaves
attached to each edge $\te$ of $\wtk_p$ and the set of leaves
attached to the edge $e=\phi(\te)$ of $\kk$ (these two sets of leaves
have the same cardinality by definition). Thus $\barkk_{\nu,p}$
can be considered as a directed $2$-complex over $\kk_\nu$.

\begin{thm}
\label{aaadg}
The kernel $\aaa$ of the natural homomorphism $\theta$ from $\dd(\kk_\nu,p)$
to $\dd(\kk,p)$ is isomorphic to the diagram group $\dd(\barkk_{\nu,p},\tp)$.
\end{thm}

\proof
Let $\phi_p\colon\dd(\wtk_p,\tp)\to\dd(\kk,p)$ and
$\bar\phi_p\colon\dd(\barkk_{\nu,p},\tp)\to\dd(\kk_\nu,p)$ be the
group homomorphisms induced by the maps $\phi$ and $\bar\phi$ (see
\cite[Section 5]{GuSa02a}). Let
$\bar\theta\colon\dd(\barkk_{\nu,p},\tp)\to\dd(\wtk_p,\tp)$ be
the homomorphism that collapses the leaves of this expansion. Then
the diagram
\be{commd}
\begin{array}{ccc}
\dd(\barkk_{\nu,p},\tp) & \stackrel{\bar\phi_p}{\longrightarrow} &
\dd(\kk_\nu,p)\vspace{0.5ex}\\
\downarrow\lefteqn{\tilde\theta} & & \downarrow\lefteqn{\theta} \\
\dd(\wtk_p,\tp)&\stackrel{\phi_p}{\longrightarrow}&\dd(\kk,p)
\end{array}
\ee is clearly commutative. The group $\dd(\wtk_p,\tp)$ is trivial
because $\wtk_p$ is a rooted $2$-tree and satisfies property T3 by
\cite[Theorem 7.2]{GuSa02a}. Hence the image of $\bar\phi_p$ is
contained in $\aaa$.

Let us show that $\bar\phi_p$ is injective. Consider any reduced
$(\tp,\tp)$-diagram $\bar\Delta$ over $\barkk_{\nu,p}$. It
suffices to prove that the diagram $\Delta=\bar\phi_p(\bar\Delta)$
has no dipoles. By contradiction, suppose that two cells $\pi_1$,
$\pi_2$ of $\Delta$ form a dipole. Then their labels are $f$,
$f^{-1}$ for some $f$ and the cells have the form $u=v$, $v=u$,
respectively, for some $u$, $v$. Let $q'$ be a path in $\Delta$
that connects $\iota(\Delta)$ and $\iota(\pi_1)=\iota(\pi_2)$.
Also let $q''$ be a path in $\Delta$ that connects
$\tau(\pi_1)=\tau(\pi_2)$ and $\tau(\Delta)$. By $r$ we denote the
path $\bott{\pi_1}=\topp{\pi_2}$ in $\Delta$. The diagram
$\bar\Delta$ has the corresponding path $\bar q'\bar r\bar q''$
from $\iota(\bar\Delta)$ to $\tau(\bar\Delta)$. The cells $\pi_1$,
$\pi_2$ have natural preimages $\bar\pi_1$, $\bar\pi_2$. We claim
that these cells form a dipole in $\bar\Delta$.

If $\pi_1$, $\pi_2$ are leaves of the form $e=e$ then their labels
are $f$ and $f^{-1}$. Hence $\pi_1$, $\pi_2$ are cells of the form
$\te=\te$ with labels $\tf$, $\tf^{-1}$. Now suppose that the
labels of $\pi_1$, $\pi_2$ belong to $\kk$. The atomic diagrams
$\ve(q')+\pi_1^{-1}+\ve(q'')$ and $\ve(q')+\pi_2+\ve(q'')$
coincide so the corresponding atomic $2$-paths (with top $1$-path
$q'rq''$) also coincide. The atomic diagrams
$\ve(\bar q')+\bar\pi_1^{-1}+\ve(\bar q'')$ and
$\ve(\bar q')+\bar\pi_2+\ve(\bar q'')$ have the same top path
$\bar q'\bar r\bar q''$. They are diagrams over $\wtk_p$ and they have
the same image under $\phi$. So they must also coincide according to U2.
Thus $\bar\Delta$ has a dipole, so it is not reduced, a contradiction.

It remains to show that $\bar\phi_p$ is surjective. Suppose that $\Delta$
is a diagram over $\kk$ that can be reduced to the trivial diagram $\ve(p)$.
Let $m$ be the number of dipoles one needs to cancel when reducing $\Delta$
to $\ve(p)$. We prove by induction on $m$ that there exists a
$(\tp,\tp)$-diagram $\wtD$ over $\wtk_p$ such that $\phi_p$ takes $\wtD$ to
$\Delta$. If $m=0$, then we have nothing to prove. Let $m\ge1$. Take the
dipole in $\Delta$ that is cancelled on the first step. Suppose it is formed
by a cell of the form $u=v$ and a cell of the form $v=u$. The diagram
$\Delta$ can be cut into two parts by a path with label of the form $q'vq''$.
Let $\Delta'$ be the result of cancelling the dipole in $\Delta$. Then
$\Delta'$ can be also cut into two parts by a path labelled by $q'uq''$. In
order to reduce this diagram, we have to cancel $m-1$ pairs of dipoles.
So we can apply the inductive assumption and find a diagram $\wtD'$ that
maps to $\Delta'$ under $\phi_p$. Consider the cut of $\wtD'$ by the path
labelled by $\tq'\tu\tq''$. Notice that there is an atomic diagram
$\ve(q')+(v=u)+\ve(q'')$ over $\kk$ whose top path is $q'vq''$. So we also
have the atomic diagram $\ve(\tq')+(\tv=\tu)+\ve(\tq'')$ over $\wtk_p$ with
top $1$-path $\tq'\tv\tq''$ (because of U2). Now we can replace the subpath
$\tu$ in $\wtD'$ by the dipole formed by the pair of mirror cells of the form
$\tu=\tv$ and $\tv=\tu$. The resulting diagram denoted by $\wtD$ clearly maps
to $\Delta$ under $\phi_p$.

Now let $\Delta_\nu$ be an arbitrary diagram in $\aaa$. By
$\Delta$ we denote its image under $\theta$ (that is, the result
of collapsing all leaves). Thus $\Delta$ represents a trivial
element of the group $\dd(\kk,p)$ and so it can be reduced to the
trivial diagram. From the previous paragraph, we know that there
exists a $(\tp,\tp)$-diagram $\wtD$ that is taken to $\Delta$ by
$\phi_p$. Passing from $\Delta_\nu$ to $\Delta$, we can collapse
one cell of the form $e=e$ at a time.

Suppose that $\Delta$ is a diagram over $\kk_\nu$ obtained from a diagram
$\wtD$ over $\barkk_{\nu,p}$ by applying $\bar\phi_p$. Also let $\Delta$
be the result of collapsing one cell of the form $e=e$ in a diagram
$\Delta_\nu$ over $\kk_\nu$. Then it suffices to show that some diagram
$\bar\Delta$ over $\barkk_{\nu,p}$ is mapped to $\Delta_\nu$ under
$\bar\phi_p$ and is mapped  to $\wtD$ under $\tilde\theta$. This is obvious:
indeed, we just find an edge $\te$ in $\wtD$ that maps to $e$, the edge that
remains after the cell $e=e$ (with some label $f$) is collapsed. If we
insert the cell $\te=\te$ with label $\tf$ into $\wtD$ instead of the edge
$\te$, then we get exactly what we need.

Thus the homomorphism $\bar\phi_p$ is in fact an isomorphism between
$\dd(\barkk_{\nu,p},\tp)$ and the subgroup $\aaa=\Ker\theta$ in
$\dd(\kk_\nu,p)$.
\endproof

Now we can use the representation of $\aaa$ as a diagram group in
order to complete the description of $\aaa$. First we show that
$\aaa$ can be represented as a diagram product. In order to do
this, we need to modify the directed $2$-complex $\barkk_{\nu,p}$
preserving its diagram groups. Let $\te$ be an edge of
$\barkk_{\nu,p}$. If $e=\phi(\te)\in\kk$ is the label of this
edge, then there are $\nu(e)$ leaves of the form $\te=\te$ in
$\barkk_{\nu,p}$. We replace these leaves by a ``switch" (see
Section \ref{dps}), that is, a $2$-cell of the form $\te=bsc$, and
$\nu(e)$ positive cells of the form $s=s$. (Here $b$, $s$, $c$
depend of $\te$.) The label of a $2$-cell $s=s$ will be the same
as the label of the leave it replaces. Applying this operation to
all edges of $\barkk_{\nu,p}$ at once, we get a new directed
$2$-complex $\hat\kk_{\nu,p}$.

\begin{lm}
\label{modif}
The diagram groups $\dd(\barkk_{\nu,p},\tp)$ and $\dd(\hat\kk_{\nu,p},\tp)$
are isomorphic.
\end{lm}

\proof
Let $\psi$ be the morphism from $\barkk_{\nu,p}$ to $\hat\kk_{\nu,p}$
that sends the leave $\te=\te$ with label $f$ to the $2$-path
$(1,\te=bsc,1)\circ(b,s=s,c)\circ(1,\te=bsc,1)^{-1}$, where the $2$-cell
$s=s$ has the same label $f$. This morphism induces a homomorphism
$\psi_{\tp}$ from $\dd(\barkk_{\nu,p},\tp)$ to $\dd(\hat\kk_{\nu,p},\tp)$
(see \cite[Section 5]{GuSa02a}). We shall show that this homomorphism is
an isomorphism.

Let us consider a reduced $(\tp,\tp)$-diagram $\Delta$ over $\hat\kk_{\nu,p}$.
Suppose that there is an edge of $\Delta$ labelled by $s$, where $\te=bsc$
is some switch of $\hat\kk_{\nu,p}$. Consider the maximal subdiagram $\Xi$
of $\Delta$ that contains this edge and is a product of $(s,s)$-cells. In
particular, $\Xi$ is an $(s,s)$-diagram. The top path of $\Xi$ is not
contained in the top path of $\Delta$ because $\topp{\Delta}$ has no edges
labelled by $s$. So there exists a cell $\pi_1$ whose bottom path contains
$\topp{\Xi}$. This cell must be a switch, that is, a $(\te,bsc)$-cell.
Analogously, we can find a $(bsc,\te)$-cell $\pi_2$ in $\Delta$ whose top
path contains $\bott{\Xi}$. We can find some paths $q_1$, $q_2$ from
$\iota(\Delta)$ to $\iota(\pi_1)$, $\iota(\pi_2)$, respectively, such that
there is a subdiagram $\Delta'$ with the top path of the form $q_1r_1$ and
the bottom path of the form $q_2r_2$. Here $r_1$ is the first edge of
$\bott{\pi_1}$ and $r_2$ is the first edge of $\topp{\pi_2}$. This implies
that $r_1$ and $r_2$ must coincide in $\Delta'$ because they are labelled by
$b$. Indeed, there are no $2$-cells whose top or bottom path ends with the
edge $b$. Therefore, the first edge of $\bott{\pi_1}$ coincides with the first
edge of $\topp{\pi_2}$. A similar argument implies that the third edge of
$\bott{\pi_1}$ coincides with the third edge of $\topp{\pi_2}$.

This shows that any edge labelled by $s$ in $\Delta$ is contained in a
$(\te,\te)$-subdiagram, which is the product of a $(\te,bsc)$-cell, some
$(s,s)$-cells, and a $(bsc,\te)$-cell. Since $\Delta$ is reduced, all
labels of the $(s,s)$-cells, read top to bottom, form a freely irreducible
word. Now we can replace our $(\te,\te)$-subdiagram by a product of
$(\te,\te)$-cells with the same labels. This operation can be done
simultaneously for all edges in $\Delta$ labelled by $s$, where $s$ is
involved in some switch of the form $\te=bsc$. The result is obviously a
preimage of $\Delta$ under $\psi_\tp$. It is easy to see that the operation
we just described defines the inverse of $\psi_\tp$, hence $\psi_\tp$ is an
isomorphism.
\endproof

It follows from Section \ref{dps} that the diagram group
$\dd(\hat\kk_{\nu,p},\tp)$ is the diagram product over $\wtk_p$ of the
family of groups $\{G_{\te}\}$, where $\te$ runs over all edges of $\wtk_p$.
Each group $G_{\te}$ is the diagram group of the directed $2$-complex
$\fff_{\nu(e)}$ with $\nu(e)$ positive $2$-cells of the form $e=e$ and base
$e$. It is clear that the connected component $\Sq(\fff_{\nu(e)},e)$ of the
Squier complex is a wedge of $\nu(e)$ circles. So $G_{\te}$ is the free
group of rank $\nu(e)$. Now we are going to prove that $\aaa$ is a partially
commutative group.

\begin{thm}
\label{aaapc}
$1$. For every edge $\te$ of $\wtk_p$, let $G_\te$ be the free group of
rank $\nu(e)$, where $e=(\phi(\te))$. Then the kernel $\aaa$ of the natural
retraction $\theta$ from $\dd(\kk_\nu,p)$ to $\dd(\kk,p)$ is a graph
product of the free groups $G_{\te}$ over the independence graph of
$\wtk_p$. Thus $\aaa$ has a partially commutative presentation, where
the generators are arbitrary symbols of the form $a(\te,i)$ with $\te$
an edge of $\wtk_p$, $1\le i\le\nu(e)$, and the defining relations are
of the form $[a(\te_1,i_1),a(\te_2,i_2)]=1$, with $\te_1$, $\te_2$
independent edges of $\wtk_p$.

$2$. As an element of the diagram group $\dd(\kk_\nu,p)$, the symbol $a(\te,i)$
is represented by a diagram of the form $\Delta^{-1}\Psi\Delta$ over $\kk_\nu$,
where $\Delta=\phi(\wtD)$ for some diagram $\wtD$ over $\wtk_p$ with
$\bott{\wtD}=\tp$ and $\topp{\wtD}=\tq'\te\tq''$ for some $\tq'$, $\tq''$, and
$\Psi$ is the atomic diagram $\ve(q')+f_i+\ve(q'')$, where $q'=\phi(\tq')$,
$q''=\phi(\tq'')$, and $f_i$ is the $i$-th leaf attached to $e=\phi(\te)$. This
diagram is independent on the choice of $\wtD$.

$3$. For every diagram $\Gamma\in\dd(\kk,p)$, and every element $a(\te,i)$
represented by the diagram $\Delta^{-1}\Psi\Delta$ as in part $2$, the
element $\Gamma^{-1}a(\te,i)\Gamma$ is equal to the generator
$a(\bar e,i)=\bar\Delta^{-1}\Psi\bar\Delta$, where $\bar\Delta$ is the
lift of the diagram $\Delta\Gamma$ over $\kk$ to the universal $2$-cover
$\wtk_p$. The edge $\bar e$ is the image of the occurrence of
$e=\phi(\te)$ in $q'eq''=\topp{\Delta\Gamma}$ after we lift the diagram
into $\wtk_p$.
\end{thm}

\proof
1. We have already shown that $\aaa$ is isomorphic to the diagram product
over $\wtk_p$ of the family of free groups $\{G_{\te}\}$ of rank $\nu(e)$.
Since $\wtk_p$ is a rooted $2$-tree (by property U1), it remains to use
Lemma \ref{dptriv}.

2. The first statement follows from properties T1 and T2, and Remark
\ref{expl}.

3. The fact that the lift of $\Delta\Gamma$ exists follows from
\cite[Lemma 8.3]{GuSa02a}. The fact that $\bar\Delta^{-1}\Psi\bar\Delta$
is a generator of $\aaa$ and the statement about $\bar e$ follows from
Remark \ref{expl}.
\endproof

This completes the description of the diagram group $\dd(\kk_\nu,p)$ as
a semi-direct product of $\aaa$ and $\dd(\kk,p)$.

\section{Representations of diagram groups by homeomorphisms}
\label{trans}

For every two natural numbers $m$, $n$, let $\bfc(m\to n)$ be the
set of all increasing homeomorphisms from the interval $[0,m]$
onto the interval $[0,n]$ of the real line. The union
$\bfc=\bigcup_{m,n}\bfc(m\to n)$ is a groupoid where the objects
are all intervals $[0,m]$ and morphisms are the functions from
$\bfc(m\to n)$. The local groups of this groupoid are the groups
$\bfc(m\to m)$.

We are also going to consider the subgroupoid $\plf$ of $\bfc$
whose morphisms are piecewise linear homeomorphisms with finitely
many singular points (breakpoints of the derivative). The set of
such piecewise linear homeomorphisms between an interval $[0,m]$
and the interval $[0,n]$ will be denoted by $\plf(m\to n)$.

Another useful subgroupoid of $\bfc$ is the groupoid $\plf_2$
where morphisms are functions from $\plf$ whose singularities
occur at dyadic rational points only and whose slopes are integer
powers of $2$. It is well known that the local groups of this
groupoid are isomorphic to R.\,Thompson's group $F$ \cite{CFP}.

In this section, we consider natural representations of diagram
groupoids of directed $2$-complexes in the groupoids $\bfc$,
$\plf$, $\plf_2$.

Let $\kk$ be any directed $2$-complex, and $\Gamma$ be one of the
groupoids $\bfc$, $\plf$, or $\plf_2$. To each $2$-cell $f$ of
$\kk$ we assign a homeomorphism $T_f$ in $\Gamma(m\to n)$, where
$m$ is the length of $\topp{f}$ and $n$ is the length of
$\bott{f}$. It is assumed that $T_{f^{-1}}=T_f^{-1}$ for every
$2$-cell $f$. We call $T_f$ a {\em transition function\/} of $f$.
The assignment of a transition function $T_f$ for each $2$-cell
$f$ is called a {\em transition scheme\/} on $\kk$ (over
$\Gamma$).

We can extend every transition scheme $T$ to the set of all
$2$-paths in $\kk$. For an atomic $2$-path $\delta=(u,f,v)$, where
$u$, $v$ are $1$-paths, $f$ is a $2$-cell, we define a function
$T_\delta$ from $[0,m]$ onto $[0,n]$, where $m=|u\topp{f}v|$,
$n=|u\bott{f}v|$ as follows:
$$
(t)T_\delta=\left\{
\begin{array}{ll}
t,&0\le t\le|u|\cr
(t-|u|)T_f,&|u|\le t\le m-|v|\cr
t+n-m,&m-|v|\le t\le m
\end{array}
\right..
$$
If $\delta=\delta_1\circ\cdots\circ\delta_r$ is a $2$-path
decomposed into a product of atomic $2$-paths, then by definition
$T_\delta$ is the composition $T_{\delta_1}\cdots T_{\delta_r}$
(it is clear that the composition is defined). This function is
called a transition function of $\delta$. It is easy to see that
this function belongs to $\Gamma(m\to n)$.

Let $\delta_1$, $\delta_2$ be independent atomic $2$-paths.
Without loss of generality, they can be written as
$\delta_1=(u,f_1,v\topp{f_2}w)$, $\delta_2=(u\topp{f_1}v,f_2,w)$
for some $2$-cells $f_1, f_2$ and $1$-paths $u, v, w$. Let
$\delta_1'=(u,f_1,v\bott{f_2}w)$,
$\delta_2'=(u\bott{f_1}v,f_2,w)$. Then the $2$-paths
$\delta_1\circ\delta_2'$ and $\delta_2\circ\delta_1'$ are
isotopic, they correspond to the same diagram
$\ve(u)+f_1+\ve(v)+f_2+\ve(w)$. It follows easily from the
definitions that the transition functions of both
$\delta_1\circ\delta_2'$ and $\delta_2\circ\delta_1'$ are equal to
the continuous function that is affine with slope $1$ on the
intervals that correspond to $u$, $v$, $w$, acts as $T_{f_1}$ on
the interval of length $|\topp{f_1}|$, and acts as $T_{f_2}$ on
the interval of length $|\topp{f_2}|$. This implies that $2$-paths
with the same diagrams have the same transition function, so any
diagram $\Delta$ over $\kk$ has a unique transition function
(depending only on the transition scheme $T$). This function has a
nice geometric description.

Suppose that $\Delta$ is a $(u,v)$-diagram over $\kk$, where
$m=|u|$, $n=|v|$. Each edge of $\Delta$ can be identified with a
closed unit interval so each point on the edge has a coordinate
from $[0,1]$. Then every $1$-path of length $r$ in $\Delta$ is
naturally identified with the interval $[0,r]$ so that each point
on this path has a coordinate between $0$ and $r$. For every cell
$\pi$ of $\Delta$ labelled by $f$, we connect each point with
coordinate $t$ on the top path of $\pi$ with the point on the
bottom path of $\pi$ that has the coordinate $(t)T_f$ on it. All
these connecting lines can be chosen disjoint since $T_f$ is
strictly increasing. If we draw the connecting lines for all the
cells of $\Delta$, then we have a disjoint ``vertical" family of
curves (a lamination) each of which consists of finitely many
connecting lines. These curves will be called the {\em transition
lines\/} of $\Delta$ (some of them may consist of just one point).
Now it is easy to describe the function $T_\Delta$: for every
point $x$ with coordinate $t\in[0,m]$ on the top path of $\Delta$,
the value of $(t)T_\Delta$ is the coordinate of the end point of
the transition line starting at $x$.

Clearly, the transition function of a dipole is the identity.
Therefore, equivalent diagrams have the same transition functions.
This implies that given a transition scheme $T$ on $\kk$ over a
groupoid $\Gamma$ (which is equal to either $\bfc$ or $\plf$ or
$\plf_2$), we have an induced functor $\psi_T$ from the diagram
groupoid $\dd(\kk)$ into the groupoid $\Gamma$: to any diagram
$\Delta$, the functor $\psi_T$ assigns the transition function
$T_\Delta$ of $\Delta$ which is an element of $\Gamma(m\to n)$,
where $m=|\topp{\Delta}|, n=|\bott{\Delta}|$.

This functor maps every diagram group $\dd(\kk,w)$ into the local
group $\Gamma(|w|\to|w|)$. Thus for every diagram group, we have a
large collection of representations by homeomorphisms of intervals
of the real line.

Some directed $2$-complexes have transition schemes which induce
faithful representations of diagram groups. This is true, for
example, in the case of the Dunce hat (see below). Notice that
since the group $\plf(m\to m)$ is totally  orderable \cite{BrSq85},
every diagram group that has enough representations in $\plf(m\to m)$
to separate all elements, is totally orderable as well. However, there
are directed $2$-complexes $\kk$ such that the intersection of kernels
of all representations induced by transition schemes of $\kk$ is not
trivial. This was mentioned in \cite[Section 17]{GuSa97} without a proof.
Now we are going to give an example.

\begin{ex}
\label{exdiag}
{\rm
Let $\kk$ be the directed $2$-complex $\la x\mid x^2=x^2\ra$ with
one vertex, one edge, and one positive $2$-cell $f$. Let $A$, $B$
be the $(x^5,x^5)$-diagrams over $\kk$ shown on Figure \theppp:

\begin{center}
\unitlength=1mm
\special{em:linewidth 0.4pt}
\linethickness{0.4pt}
\begin{picture}(140.00,30.00)
\put(0.00,17.00){\vector(1,0){15.00}}
\put(15.00,17.00){\vector(1,1){10.00}}
\put(25.00,27.00){\vector(1,0){15.00}}
\put(40.00,27.00){\vector(1,0){15.00}}
\put(55.00,27.00){\vector(1,-1){10.00}}
\put(15.00,17.00){\vector(1,0){15.00}}
\put(30.00,17.00){\vector(1,0){15.00}}
\put(45.00,17.00){\vector(1,1){10.00}}
\put(45.00,17.00){\vector(1,-1){10.00}}
\put(55.00,7.00){\vector(1,1){10.00}}
\put(30.00,17.00){\vector(1,-1){10.00}}
\put(30.00,17.00){\vector(1,1){10.00}}
\put(40.00,7.00){\vector(1,0){15.00}}
\put(15.00,17.00){\vector(1,-1){10.00}}
\put(25.00,7.00){\vector(1,0){15.00}}
\put(75.00,17.00){\vector(1,1){10.00}}
\put(85.00,27.00){\vector(1,-1){10.00}}
\put(95.00,17.00){\vector(1,0){15.00}}
\put(110.00,17.00){\vector(1,0){15.00}}
\put(125.00,17.00){\vector(1,0){15.00}}
\put(75.00,17.00){\vector(1,-1){10.00}}
\put(85.00,7.00){\vector(1,1){10.00}}
\put(27.00,22.00){\makebox(0,0)[cc]{$f$}}
\put(42.00,22.00){\makebox(0,0)[cc]{$f^{-1}$}}
\put(55.00,17.00){\makebox(0,0)[cc]{$f$}}
\put(85.00,17.00){\makebox(0,0)[cc]{$f$}}
\put(27.00,12.00){\makebox(0,0)[cc]{$f^{-1}$}}
\put(42.00,12.00){\makebox(0,0)[cc]{$f$}}
\put(40.00,1.00){\makebox(0,0)[cc]{\Large $A$}}
\put(114.00,1.00){\makebox(0,0)[cc]{\Large $B$}}
\put(20.00,18.00){\makebox(0,0)[cc]{\tiny1}}
\put(35.00,18.00){\makebox(0,0)[cc]{\tiny2}}
\put(48.00,17.00){\makebox(0,0)[cc]{\tiny3}}
\put(35.00,15.00){\makebox(0,0)[cc]{\tiny4}}
\put(20.00,15.00){\makebox(0,0)[cc]{\tiny5}}
\end{picture}
\end{center}
\begin{center}
\nopagebreak[4] Figure \theppp.
\end{center}
\addtocounter{ppp}{1}
\noindent
Thus $A$ corresponds to the $2$-path $(x,f,x^2)\circ(x^2,f^{-1},x)
\circ(x^3,f,1)\circ(x^2,f,x)\circ(x,f^{-1},x^2)$ whereas $B$
corresponds to the $2$-path $(1,f,x^3)$. The cells of $A$ are
enumerated according to their appearance in the $2$-path. Let
$T_f$ be any strictly increasing continuous function from $[0,2]$
onto itself. First we show that the transition function $T_A$ from
$[0,5]$ onto itself is identical on $[0,2]$. If $t\in[0,1]$, then
$(t)T_A=t$. So let $t\in[1,2]$. Suppose that $(t-1)T_f\le1$. In
this case the transition line of the point with coordinate $t$ on
the top of $A$ will consist of two paths that go through the first
and the fifth cell of $A$. These paths are mirror images of each
other. So $(t)T_A=1+(t-1)T_fT_f^{-1}=t$ in this case. Now suppose
that $(t-1)T_f>1$. This means that the part of the transition line
of $t$ inside the first cell of $A$ ends on some point on the top
path of the second cell of $A$ and its coordinate there will be
some number $s=(t-1)T_f-1\le1$. The portion of the transition line
inside the second cell of $A$ will connect the point with
coordinate $s$ on the top to the point with coordinate $(s)T_f^{-1}$
on the bottom. Since $T_f$ is increasing,
$(s)T_f^{-1}\le(1)T_f^{-1}<t-1\le1$ so the line goes through the
common boundary of the second and the fourth cell. Therefore, it
also consists of two parts that are mirror images of each other,
and so $(t)T_A=t$.

We proved that $T_A$ is identical on $[0,2]$. It is obvious that
$T_B$ is identical on $[2,5]$. Thus $T_A$ and $T_B$ commute. Hence
the transition function of the commutator $C=[A,B]$ is identical
on the whole interval $[0,5]$. The diagram $C$ has $12$ cells and
it is clearly reduced. So $C$ is an example of a diagram that is
contained in the kernel of any homomorphism induced by a
transition scheme.
}
\end{ex}

Notice that the diagram group $G$ of the complex $\la x\mid x^2=x^2\ra$
with base $x^5$ is a partially commutative group given by
$\la a,b,c,d\mid[a,b]=[b,c]=[c,d]=1\ra$. The component $\Sq(\kk,x^5)$
has only one vertex $x^5$ and four positive edges $a=(x,f,x^2)$,
$b=(x^3,f,1)$, $c=(1,f,x^3)$, $d=(x^2,f,x)$. There are three $2$-cells
(squares $\ve(x)+f+f+\ve(1)$, $\ve(1)+f+\ve(x)+f+\ve(1)$,
$\ve(1)+f+f+\ve(x)$) which correspond to the three commutativity relations.

It is obvious that the groups $\bfc(n\to n)$ (resp., $\plf(n\to n)$)
are isomorphic to each other for all $n>0$. A more traditional notation
for the groups $\bfc(1\to1)$ and $\plf(1\to1)$ is $\Homeo_+[0,1]$ and
$\plf_+[0,1]$, respectively. Example \ref{exdiag} shows that homomorphisms
into $\Homeo_+[0,1]$ induced by transition schemes do not always separate
elements of a diagram group. Thus it is natural to ask the following
question (a similar question was mentioned in \cite{GuSa97}).

\begin{prob}
\label{pr3}
Is it true that any diagram group is residually $\Homeo_+[0,1]$ or even
residually $\PLF_+[0,1]$?
\end{prob}

Notice that we also do not know if any diagram group is residually $F$.
\vspace{1ex}

Consider the Dunce hat $\hh_0=\la x\mid x=x^2\ra$. In this case
every transition scheme consists of one function $h\colon[0,1]\to[0,2]$.
Let $h$ be the function $t\mapsto 2t$, and $T$ be the corresponding
transition scheme. We shall show, in particular, that the induced
homomorphism of the R.\,Thompson group $F=\dd(\hh_0,x)$ into $\plf_2$
is the well known representation of $F$ by piecewise linear transformations
of the unit interval \cite{CFP}.

\begin{lm}
\label{linhom}
The homomorphism $\psi_T$ is injective on $\dd(\hh_0,x)$.
\end{lm}

\proof Consider the following two diagrams over $\hh_0$:

\begin{center}
\unitlength=.55mm \special{em:linewidth 0.4pt}
\linethickness{0.4pt}
\begin{picture}(94.00,40.00)
\put(3.00,24.00){\circle*{1.33}}
\put(13.00,24.00){\circle*{1.33}}
\put(23.00,24.00){\circle*{1.33}}
\put(33.00,24.00){\circle*{1.33}}
\put(53.00,24.00){\circle*{1.33}}
\put(63.00,24.00){\circle*{1.33}}
\put(73.00,24.00){\circle*{1.33}}
\put(83.00,24.00){\circle*{1.33}}
\put(93.00,24.00){\circle*{1.33}}
\put(3.00,24.00){\line(1,0){10.00}}
\put(13.00,24.00){\line(1,0){10.00}}
\put(23.00,24.00){\line(1,0){10.00}}
\put(53.00,24.00){\line(1,0){10.00}}
\put(63.00,24.00){\line(1,0){10.00}}
\put(73.00,24.00){\line(1,0){10.00}}
\put(83.00,24.00){\line(1,0){10.00}}
\bezier{120}(13.00,24.00)(23.00,35.00)(33.00,24.00)
\bezier{120}(3.00,24.00)(12.00,13.00)(23.00,24.00)
\bezier{256}(3.00,24.00)(19.00,52.00)(33.00,24.00)
\bezier{256}(3.00,24.00)(14.00,-4.00)(33.00,24.00)
\bezier{132}(63.00,24.00)(71.00,37.00)(83.00,24.00)
\bezier{108}(53.00,24.00)(64.00,15.00)(73.00,24.00)
\bezier{208}(53.00,24.00)(70.00,45.00)(83.00,24.00)
\bezier{176}(53.00,24.00)(65.00,8.00)(83.00,24.00)
\bezier{296}(53.00,24.00)(67.00,55.00)(93.00,24.00)
\bezier{296}(53.00,24.00)(62.00,-6.00)(93.00,24.00)
\put(18.00,2.00){\makebox(0,0)[cc]{{\Large $x_0$}}}
\put(73.00,2.00){\makebox(0,0)[cc]{{\Large $x_1$}}}
\end{picture}
\nopagebreak[4]

Figure \theppp.
\end{center}
\addtocounter{ppp}{1}

It was mentioned in \cite[Section 1]{GuSa99} that these diagrams
generate the group $\dd(\hh_0,x)$.

The straightforward computation gives
$(1/2)T_{x_0}T_{x_1}=(1/4)T_{x_1}=1/8$ and
$(1/2)T_{x_1}T_{x_0}=(1/2)T_{x_0}=1/4$. So the transition
functions $T_{x_0}$ and $T_{x_1}$ do not commute. This means that
the image of $\dd(\hh_0,x)$ under $\psi_T$ is non-Abelian. It is
well known, however, that all proper homomorphic images of the
group $F$ are Abelian \cite{CFP}. Thus $\psi_T$ is injective on
$\dd(\hh_0,x)$.
\endproof

\begin{lm}
\label{convdiag} The functor $\psi_T$ is an isomorphism from
$\dd(\hh_0)$ onto $\plf_2$.
\end{lm}

\proof
Let $h$ be a function from $\plf_2(m\to n)$. We need to show
that there exists a unique reduced $(x^m,x^n)$-diagram $\Delta$
over $\hh_0$ such that the transition function $T_\Delta$
coincides with $h$.

We begin with proving the uniqueness. Let $\Delta_1$, $\Delta_2$
have the same transition function $h$. Since each equivalence
class of diagrams has exactly one reduced representative
\cite[Theorem 2.5]{GuSa02a}, it suffices to prove that $\Delta_1$
and $\Delta_2$ are equivalent. Both diagrams are
$(x^m,x^n)$-diagrams so we can consider the $(x^n,x^n)$-diagram
$\Delta_1^{-1}\Delta_2$. The transition function of it is the
identity. Hence the same is true for the transition function of
$\Theta=\Delta^{-1}\Delta_1^{-1}\Delta_2\Delta$, where $\Delta$ is
any $(x^n,x)$-diagram over $\hh_0$. By Lemma \ref{linhom}, the
diagram $\Theta$ represents the trivial element of $\dd(\hh_0,x)$
so it is equivalent to $\ve(x)$. It follows immediately that
$\Delta_1$, $\Delta_2$ are equivalent.

Now let us prove the existence. Given a function $h$ from
$\plf_2(m\to n)$, let $a_0=0$, $a_1$, \dots, $a_k$, $a_{k+1}=m$ be
the increasing sequence that contains all the singularity points
of $h$ ($k\ge0$). Let $b_i=(a_i)h$ for all $0\le i\le k+1$. By $d$
we denote a smallest natural number such that all the numbers
$2^da_i$, $2^db_i$ are integers ($0\le i\le k+1$). Consider the
function $g$ from $[0,2^dm]$ onto $[0,2^dn]$ given by
$(t)g=2^d\cdot(2^{-d}t)h$. All singularities of $g$ occur at
integer points and all slopes are the same as the corresponding
slopes of $h$. The map $g$ takes $[2^da_i,2^da_{i+1}]$ onto
$[2^db_i,2^nd_{i+1}]$ for all $0\le i\le k$. On each of these
intervals, $g$ has some slope of the form $2^{c_i}$, where $c_i$
is an integer.

For any $j\ge0$, let $\Delta_j$ be an $(x,x^{2^j})$-diagram over
$\hh_0$ that can be defined by induction as follows. We let
$\Delta_0=\ve(x)$ and $\Delta_{j+1}=\pi\circ(\Delta_j+\Delta_j)$
for all $j\ge0$. Here $\pi$ is the $(x,x^2)$-diagram of one cell.
The transition function of $\Delta_j$ is obviously linear. Now for
any number $r\ge1$, we may take the sum of $r$ copies of $\Delta_j$.
We will denote this sum by $r\cdot\Delta_j$. It corresponds to the
linear function from $[0,r]$ onto $[0,2^jr]$.

For any $0\le i\le k$, if $d_i\ge0$, then we let $r_i=2^n(a_{i+1}-a_i)$.
If $d_i<0$, then we let $r_i=2^n(b_{i+1}-b_i)$. By definition, the
diagram $\Xi$ is the sum of diagrams $\Xi_i$ ($0\le i\le k$), where
$\Xi_i=r_i\cdot\Delta_{d_i}$, if $d_i\ge0$ and
$\Xi_i=(r_i\cdot\Delta_{-d_i})^{-1}$, if $d_i<0$. The transition
function of $\Xi$ is exactly $g$. Now by $\Delta$ we will denote
the reduced form of the diagram
$(m\cdot\Delta_d)\circ\Xi\circ(n\cdot\Delta_d)^{-1}$, whose
transition function is $h$.
\endproof

\begin{rk}
\label{cinfty} {\rm One can also consider representations of
diagram groups by differentiable functions. For example, consider
subgroupoid of $\bfc$ that consists of $
C^\infty$-homeo\-morph\-isms. It is easy to see that if every
function in a transition scheme of a directed $2$-complex $\kk$ is
$C^\infty$ and has derivative $1$ at the ends of its domain then
the induced representation maps the diagram groupoid into groupoid
of $C^\infty$-homeomorphisms. In particular, one can construct
many faithful $C^\infty$-representations of the R.\,Thompson group
$F$. The existence of such representations was first proven in
\cite{GhSe}. }
\end{rk}

\begin{rk}
\label{rkhom} {\rm Notice that one can consider transition schemes
of directed $2$-complexes over other groupoids. The groupoid
should only have an operation of tensor product \cite{JS}. The
groupoid $\bfc$ has a natural tensor product: if $f\in \bfc(k\to
l)$, $g\in \bfc(m\to n)$ then the tensor product of $f$ and $g$ is
the function from $\bfc(k+m\to l+m)$ which acts on $[0,k]$ as $f$
and on $[k,k+m]$ as $g(x-k)+l$. Another example of a groupoid with
a tensor product is, of course, the diagram groupoid of any
directed $2$-complex with one vertex. In that case the tensor
product is the sum of diagrams.

Instead of $\bfc$, one can consider the groupoid of all
homeomorphisms between rectangles of the form $[0,m]\times[0,1]$
on the plane $\mathbb{R}^2$ ($m\in \mathbb{N}$), that map isometrically $\{0\}\times [0,1]$ to $\{0\}\times[0,1]$ and $\{m\}\times [0,1]$ to $\{n\}\times[0,1]$. This groupoid also has a natural operation of
tensor product. In that case a transition scheme $T$ of a directed
$2$-complex $\kk$ associates with every $2$-cell $f$ a
homeomorphism from $[0,m]\times[0,1]$ onto $[0,n]\times[0,1]$,
where $m=|\topp{f}|$, $n=|\bott{f}|$ . We can prove that the
functors induced by these transition schemes separate all elements
of any diagram groupoid. }
\end{rk}

\section{The structure of $\gr_1$}
\label{strg-1}

The complex $\hh_1$ is an expansion of the Dunce hat
$\hh_0=\la x\mid x^2=x\ra$ with $\nu(x)=1$. Hence we can apply the results
of Section \ref{expan}. In this section, we are going to give a precise
description of the group $\gr_1=\dd(\hh_1,x)$.

As in Section \ref{expan}, let $\theta$ be the (natural) retraction
from $\gr_1=\dd(\hh_1,x)$ to $F\cong\dd(\hh_0,x)$, the group $\aaa$ be
the kernel of this retraction, $\wtk=\wtk_x$ be the universal $2$-cover
of $\kk=\hh_0$ with base $x$. Our first goal is to give a nice description
of $\wtk$. We are going to construct a certain directed $2$-complex $\mm$
over $\hh_0$, and then show that it satisfies properties U1 and U2$'$
for $\kk=\hh_0$. Notice that we do not need to describe the labelling
morphism $\phi$ because $\hh_0$ has only one edge and only one positive
$2$-cell.

The set of vertices of $\mm$ is the set of all dyadic rational numbers in
the unit interval $[0;1]$. The complex $\mm$ is a union of complexes
$\mm_i$ ($i\ge 0$). All these complexes have the same set of vertices.
(Warning: this is not a natural filtration of $\mm$.)

The complex $\mm_0$ is defined by the following process. We start with an
edge $\tx$ that goes from $0$ to $1$. This is the edge of {level\/} $0$.
Suppose that all edges of level $r\ge0$ have been constructed. To each edge
$e$ of level $r$ we assign two new edges, $e'$ and $e''$. If $e$ goes from
$\alpha$ to $\beta$, where $0\le\alpha<\beta\le1$ are dyadic rational numbers,
then $e'$ goes from $\alpha$ to $(\alpha+\beta)/2$ and $e''$ goes from
$(\alpha+\beta)/2$ to $\beta$. We also attach a positive $2$-cell $f$ of the
form $e=e'e''$. By definition, we say that the $2$-cell $f$ and the edges
$e'$, $e''$ have level $r+1$.

This inductive process creates $2^r$ edges for each level $r\ge0$ and the
same number of positive $2$-cells that have level $r+1$. This gives us a
complex $\mm_0$ that consists of all edges and $2$-cells constructed during
this process. We can draw a part of $\mm_0$ in the following picture.

\begin{center}
\unitlength=0.7mm
\special{em:linewidth 0.4pt}
\linethickness{0.4pt}
\begin{picture}(104.00,48.00)
\put(8.00,35.00){\circle*{1.33}}
\put(88.00,35.00){\circle*{1.33}}
\bezier{348}(8.00,35.00)(48.00,52.00)(88.00,35.00)
\bezier{400}(8.00,35.00)(48.00,5.00)(88.00,35.00)
\put(48.00,20.00){\circle*{1.33}}
\bezier{316}(8.00,35.00)(17.00,-4.00)(48.00,20.00)
\bezier{336}(48.00,20.00)(84.00,-5.00)(88.00,35.00)
\put(21.00,12.00){\circle*{1.33}}
\put(76.00,11.00){\circle*{1.33}}
\bezier{248}(8.00,35.00)(-10.00,10.00)(21.00,12.00)
\bezier{232}(21.00,12.00)(38.00,-9.00)(48.00,20.00)
\bezier{236}(48.00,20.00)(59.00,-10.00)(76.00,11.00)
\bezier{220}(76.00,11.00)(104.00,14.00)(88.00,35.00)
\put(2.00,18.00){\circle*{1.33}}
\put(35.00,3.00){\circle*{1.33}}
\put(60.00,3.00){\circle*{1.33}}
\put(94.00,19.00){\circle*{1.33}}
\put(48.00,47.00){\makebox(0,0)[cc]{\Large $\tx$}}
\put(7.00,39.00){\makebox(0,0)[cc]{$0$}}
\put(91.00,38.00){\makebox(0,0)[cc]{$1$}}
\put(48.00,23.00){\makebox(0,0)[cc]{$^{1/2}$}}
\put(23.00,15.00){\makebox(0,0)[cc]{$^{1/4}$}}
\put(76.00,14.00){\makebox(0,0)[cc]{$^{3/4}$}}
\put(7.00,18.00){\makebox(0,0)[cc]{$^{1/8}$}}
\put(35.00,6.00){\makebox(0,0)[cc]{$^{3/8}$}}
\put(62.00,6.00){\makebox(0,0)[cc]{$^{5/8}$}}
\put(90.00,20.00){\makebox(0,0)[cc]{$^{7/8}$}}
\end{picture}

\nopagebreak[4] Figure \theppp.
\end{center}
\addtocounter{ppp}{1}

\noindent
It shows all $2$-cells and all edges of level $\le3$. It is
convenient to view $\mm_0$ as a plane complex. By definition, we
set the {\em height\/} of all edges and $2$-cells of $\mm_0$ to
$0$. For any $h\ge1$, assuming that $\mm_{h-1}$ is already
constructed, we consider all pairs of consecutive edges $e'$, $e''$
in $\mm_{h-1}$ (this means that $e'e''$ is a $1$-path). For each
of these pairs, if $\mm_{h-1}$ contains no $2$-cells of the form
$e=e'e''$ , where $e$ is an edge, we add a new edge $e$ and a new
$2$-cell of the form $e=e'e''$. We set the height of all edges and
$2$-cells added in this way as $h$. The new directed $2$-complex is
denoted by $\mm_h$.

By definition, $\mm=\bigcup_{h\ge0}\mm_h$ is a rooted $2$-tree with
root $\tx$ so property U1 holds.

Let $\E$ be the set of edges of $\mm$. It is easy to see that for each edge
$e\in\E$, there is a unique $2$-cell of the form $e=e'e''$, where
$e',e''\in\E$. This cell depends of $e$ only so we denote it by $\pi_e$.
Thus given an edge $e$, we have two functions $e\mapsto e'$ and
$e\mapsto e''$ from from $\E$ to itself. Moreover, for any pair of
consecutive edges of $\mm$ we have a unique edge $e$ in $\mm$ such
that these consecutive edges are exactly $e'$, $e''$. These properties
can also be used to define $\mm$ axiomatically. They immediately imply
that $\mm$ satisfies U2$'$. So $\mm$ is the universal $2$-cover of the Dunce
hat $\hh_0$ with base $x$ according to \cite[Remark 8.2]{GuSa02a}. Notice
also that one can regard the edges of $\mm$ as elements of the one-generated
free Cantor algebra (see \cite{Hig}).

From now on, we will denote the universal $2$-cover of the Dunce hat with
base $x$ by $\wtk$. As in the Section \ref{expan}, $\barkk=\barkk_{\nu,x}$
is an expansion of $\wtk$ obtained by adding a leaf to each edge of $\wtk$.
Theorem \ref{aaadg} implies:

\begin{cy}
\label{aa-dg}
The kernel $\aaa$ of the retraction $\theta$ of $\gr_1$ onto $F$ is the
diagram group of the directed complex $\barkk$ with base $\tx$.
\end{cy}

Theorem \ref{aaapc} gives the following partially commutative presentation
of $\aaa$. We will denote the generators of $\aaa$ by $a(z)$ instead of
$a(z,i)$ as in Theorem \ref{aaapc} because in this case $i=1$ always. By
definition of $\wtk$, the endpoints of every  edge $z$ in $\wtk$ are dyadic
rational numbers with $0\le\iota(e)<\tau(e)\le1$, and for every two dyadic
numbers $0\le\lambda<\mu\le1$, there is an edge connecting $\lambda$ with $\mu$
(in fact there are countably many such edges of different heights). Thus
every edge in $\E$ corresponds to an open interval $(\iota(e),\tau(e))$ on
the real line. We claim that (by Theorem \ref{aaapc}) symbols $a(z_1)$,
$a(z_2)$ commute if and only if the corresponding intervals do not
intersect. Indeed, if $z_1$ and $z_2$ occur in the same $1$-path $q$ in
$\wtk$ from $0=\iota(\tx)$ to $1=\tau(\tx)$, then we can assume without
loss of generality that $z_1$ is to the left of $z_2$ in that path. This
means $\tau(z_1)\le \iota(z_2)$ and the intervals do not intersect.
Conversely, if the two open intervals do not intersect, then without loss
of generality $\tau(z_1)\le \iota(z_2)$. Then there exists a path in $\wtk$
(of length at most $5$) that connects vertices $0$ and $1$ and contains
$z_1$ and $z_2$, so the symbols $a(z_1)$, $a(z_2)$ commute.

So we have the following structural description of the subgroup $\aaa$ of
$\gr_1$.

\begin{thm}
\label{aaadscr}
For every open subinterval $P$ of $(0,1)$ with dyadic rational endpoints
let $Z(P)$ be a countable set of symbols. Then $\aaa$ is generated by the
union of all $Z(P)$ subject to commutativity relations: two symbols commute
if and only in the corresponding intervals do not intersect.
\end{thm}

Since $\gr_1$ is a semi-direct product of $\aaa$ and
R.\,Thompson's group $F$, it remains to describe the action of $F$
on $\aaa$. We know that the generators of $\aaa$ are in one-to-one
correspondence with edges of $\wtk$. Now to every edge of $\wtk$,
we assign a piecewise linear function. In order to do that, we
will use transition functions defined in the previous section. We
shall always assume that the transition scheme $T$ on $\hh_0$
consists of linear functions.

We need some useful technical lemma. Recall that $\phi$ denotes the
labelling function from $\wtk$ to $\hh_0$. Applying $\phi$ to a diagram,
means that we replace all its edge labels by $x$. (One can forget about
the labels of cells because they can be easily recovered.) We also recall
that vertices of $\wtk$ are dyadic rational numbers from $[0,1]$.

\begin{lm}
\label{images} $1$. Suppose that $\Delta$ is a diagram over
$\hh_0$ that has cells of the form $x^2=x$ only. Then the
transition function $T_\Delta$ is linear on each edge of the top
path of $\Delta$.

$2$. Let $\wtD$ be a $(\tq,\tx)$-diagram over $\wtk$. Suppose that $\tilde o$
is a vertex on $\tq$. Let $\Delta=\phi(\wtD)$ and let $o$ be the image
of $\tilde o$ in $\Delta$. Then $o$ is connected by a transition line in
$\Delta$ to the point on $\bott{\Delta}$ that has coordinate $\tilde o$.
\end{lm}

\proof To prove part 1, we proceed by induction on the number of
cells in $\Delta$. If $\Delta$ has no cells, then the transition
function of it is identical. Suppose that $\Delta$ has cells. In
this case there is a cell $\pi$ in $\Delta$ such that $\topp{\pi}$
is contained in $\topp{\Delta}$. Let $\Delta'$ be a subdiagram in
$\wtD$ that is obtained from $\Delta$ by deleting $\pi$ together
with $\topp{\pi}$. If the edge $e$ on the top of $\Delta$ is not
contained in $\topp{\pi}$, then the transition functions of
$\Delta$ and $\Delta'$ coincide on the unit interval $e$. Hence
the result follows by the inductive assumption. Now suppose that
the edge $e$ is contained in $\topp{\pi}$. Since the transition
function of $\pi$ is linear, the top of $\pi$ maps $e$ linearly to
a half of the edge $e'=\bott{\pi}$. The transition function of
$\Delta'$ is linear on $e'$ by the inductive assumption. It
remains to compose two linear functions.

To prove part 2, we can first assume that $\wtD$ is reduced. By
the properties of $\wtk$, the diagram $\Delta$ is also reduced. So
according to \cite[Example 2]{GuSa99}, it can be decomposed as
$\Delta_{+}\circ\Delta_{-}$, where $\Delta_{+}$ has only cells
labelled by $x=x^2$ and $\Delta_{-}$ has only cells labelled by
$x^2=x$. In addition, all vertices of $\Delta$ belong to the path
$r=\bott{\Delta_{+}}=\topp{\Delta_{-}}$. The transition line of
$o$ is thus contained in $\Delta_{-}$.

Since $\phi(\wtD)=\Delta$, the diagram $\wtD$ has an induced
decomposition $\wtD=\wtD_{+}\circ\wtD_{-}$ and the point $\tilde
o$ belongs to the path $\tr$ that cuts $\wtD$ into the two
subdiagrams. For any cell of $\wtD_{-}$, its top path has length
$2$ and its bottom path has length $1$. Since the bottom of
$\wtD_{-}$ is $\tx$, it is easy to see by induction on the number
of cells that $\wtD'$ is a diagram over the directed $2$-complex
$\mm_0$. One can lift the transition line of $o$ in $\Delta_{-}$
into $\wtD_{-}$. We just need to show that this line connects
$\tilde o$ with the point on $\tx$ that has the coordinate $\tilde
o$. This is obvious if $\wtD_{-}$ has no cells. Let us proceed by
induction on the number of cells in $\wtD_{-}$. Let $\pi$ be a
cell in $\wtD_{-}$ such that $\topp{\pi}$ is contained in
$\topp{\wtD_{-}}$. Let $\wtD'$ be a subdiagram in $\wtD_{-}$
obtained by deleting $\pi$ together with $\topp{\pi}$.

If $\tilde o$ is not an inner point of $\topp{\pi}$, then $\tilde
o$ is contained in $\wtD'$. Hence the inductive assumption can be
applied to $\wtD'$.

Now let $\tilde o$ be an inner point of $\topp{\pi}$. In this case
$\topp{\pi}$ is a path of length $2$ of the form $\te_1\te_2$,
where $\te_1$ connects some vertex $\tilde o_1$ with $\tilde o$
and $\te_2$ connects $\tilde o$ with some vertex $\tilde o_2$. The
inductive assumption applied to $\wtD'$ allows us to conclude that
the points $\tilde o_j$ are connected by their transition lines in
$\wtD'$ to the points on $\tx=\bott{\Delta'}=\bott{\Delta}$ with
coordinates $\tilde o_j$ ($j=1,2$). By definition of the complex
$\mm_0$, one has $\tilde o=(\tilde o_1+\tilde o_2)/2$. The
transition line of $\tilde o$ connects it first to the midpoint of
the edge $\bott{\pi}$. Applying part 1 to $\wtD'$, we see that the
transition function of $\wtD'$ has to be linear on the edge
$\bott{\pi}$. So the midpoint of this edge will be taken to the
midpoint of the interval between the images of $\tilde o_1$ and
$\tilde o_2$ on $\tx$. By the inductive assumption, the latter
points have coordinates $\tilde o_1$, $\tilde o_2$. Thus the midpoint
is $\tilde o$, as desired.
\endproof

Let $\te$ be an edge of $\wtk$. First we choose any $1$-paths
$\tq'$, $\tq''$ in $\wtk$ such that the $1$-path $\tq'\te\tq''$
connects $0$ and $1$. By Properties T2 and T3, there exists a
unique reduced $(\tq'\te\tq'',\tx)$-diagram over $\wtk$. If we
apply $\phi$ to it, then we get a diagram $\Delta$ over $\hh_0$.
The top path of $\Delta$ is decomposed as $q'eq''$ in a natural
way, that is, $e$ is the image of $\te$. Let $e$ be the $k$th edge
of $\topp{\Delta}$ and let $n=|q'eq''|$. The points with
coordinates $k-1$ and $k$ (the endpoints of $e$) are taken to some
dyadic rationals $\lambda$ and $\mu$ via the transition function
of $\Delta$. It follows from Lemma \ref{images} that these numbers
are exactly the endpoints of $\te$ in $\wtk$. Consider the
function $g\colon[0,1]\to [0,1]$ defined by
$(t)g=(t+k-1)T_\Delta$. We associate this function with the edge
$\te$. Let us prove that $g$ does not depend on the choice of
$\tq'$ and $\tq''$.

Suppose that we change the paths $\tq'$, $\tq''$. Then $\Delta$
also changes to some diagram, which will be equivalent to the
product of some diagram $\Psi=\Psi_1+\ve(x)+\Psi_2$ and $\Delta$,
where $\Psi_1$, $\Psi_2$ are spherical. The transition function of
$\Psi$ takes $[0,k-1]$ and $[k,n]$ to themselves and it is
identical on $[k-1,k]$. So the restriction of $T_\Delta\Psi$ on
$[k-1,k]$ will coincide to the restriction of $T_\Delta$ on
$[k-1,k]$. Thus we can say that to any edge $\te$ of $\wtk$ we
assign a unique function $g_{\te}$ from $[0,1]$ to itself. By
Lemma \ref{convdiag}, this is a piecewise linear homeomorphism
from $[0,1]$ to $[\lambda,\mu]$ with dyadic $\lambda$, $\mu$,
finitely many dyadic singularities and all slopes of the form
$2^k$, $k\in\mathbb{Z}$.

Let us denote by $\Phi$ the set of all such functions from $[0,1]$
to itself. It is easy to see that $\Phi$ is a semigroup under
composition. We are going to show that for any $h\in\Phi$, there
exists a unique edge $\te$ of $\wtk$ such that $h=g_{\te}$.

\begin{lm}
\label{1t1}
There is a natural one-to-one correspondence between edges of $\wtk$
and continuous piecewise linear functions from $[0,1]$ onto a
subinterval in $[0,1]$ with dyadic rational endpoints. All these
functions have singularities at finitely many dyadic rational points
only and all slopes of these functions are integer powers of $2$.
\end{lm}

\proof
Let $h\in\Phi$. We can write $(0)h=l/2^d$, $(1)h=m/2^d$ for some
integers $l$, $m$, $d$. Let us consider the function
$$
(t)\bar h=\left\{
\begin{array}{ll}
t/2^d,&0\le t\le l\cr
(t-l)h,&l\le t\le l+1\cr
1+(t-n)/2^d,&l+1\le t\le n
\end{array}
\right.
$$
from $[0,n]$ onto $[0,1]$, where $n=2^d+l+1-m$. It clearly belongs
to $\plf_2(n\to1)$. By Lemma \ref{convdiag}, there exists a
reduced $(x^n,x)$-diagram $\Delta$ over $\hh_0$ that has $\bar h$
as its transition function. Using \cite[Lemma 8.3]{GuSa02a}, we
can lift $\Delta$ to $\wtk$. This gives us a $(\tq,\tx)$-diagram
$\wtD$. Since $m\le2^d$, one has $|\tq|=n\ge l+1$. So let $\te$ be
the $(l+1)$th edge of $\tq$. Then $g_{\te}=h$. Indeed, if we
replace all edges of $\wtD$ by $x$, we get the diagram $\Delta$ by
the definition of the lift. Its transition function is $\bar h$,
the restriction of this function to the unit interval $[l,l+1]$ is
the function $h$ with argument shifted by $l$.

It remains to prove the uniqueness. Suppose that $\te$, $\te'$ are
two edges of $\wtk$ such that $g_{\te}=g_{\te'}=h$. By Lemma \ref{images},
both edges have the same endpoints $(0)h$, $(1)h$. So they are homotopic
in $\wtk$ by property T2. One can choose a reduced $(\te',\te)$-diagram
$\tilde\Gamma$ over $\wtk$. Let us choose any $1$-path $\tq$ in $\wtk$
from $0$ to $1$ that contains $\te$ and some $(\tq,\tx)$-diagram $\wtD$
over $\kk$. Then one can take the diagram
$\wtD'=(\tilde\Psi_1+\tilde\Gamma+\tilde\Psi_2)\wtD$ over $\wtk$
whose top path contains $\te'$, where $\tilde\Psi_1$,
$\tilde\Psi_2$ have no cells. Now we replace all edge labels in
$\wtD$, $\wtD'$, $\tilde\Gamma$, $\tilde\Psi_j$ by $x$ and get
diagrams $\Delta$, $\Delta'$, $\Gamma$, $\Psi_j$ ($j=1,2$),
respectively. Notice that $\Gamma$ is reduced because so was
$\tilde\Gamma$. The transition function of $\Delta'$ is the
composition of the transition function of $\Psi_1+\Gamma+\Psi_2$
and $\Delta$. If we restrict these functions on the corresponding
unit interval, then we see that $g_{\te'}$ will be a composition
of $g$ and $g_{\te}$, where $g$ is the transition function of
$\Gamma$. Thus $g$ must be identical so $\Gamma$ has no cells by
Lemma \ref{linhom}. But this means that $\te$, $\te'$ coincide.
\endproof

Now it is natural to change the notation for the generators of
$\aaa$. We replace the symbol $a(\te)$ by the symbol $\alpha_h$,
where $h=g_{\te}$. So $\aaa$ is generated by the symbols of the
form $\alpha_h$, where $h$ runs over $\Phi$. According to Theorem
\ref{aaadscr}, two generators $\alpha_{h_1}$, $\alpha_{h_2}$
commute if and only if the images of $h_1$ and $h_2$ have disjoint
interiors.

Now we are ready to describe the action of $F$ on $\aaa$. We can
think about the $\alpha_h$'s as of elements in $\gr_1$. The group
$F$, as a subgroup of $\gr_1$, acts on $\aaa$ by conjugation.
Recall that $F$ is isomorphic to $\plf_2(1\to1)$ that is a
subgroup of the semigroup $\Phi$. So we shall view $F$ as a
subgroup of $\Phi$. The action of $F$ is described in the
following statement.

\begin{lm}
\label{act} For any $h\in\Phi$ and for any $g\in F$, one has $g^{-1}
\alpha_hg=\alpha_{hg}$ in the group $\gr_1$.
\end{lm}

\proof First we need to find out what are the diagrams over
$\hh_1$ that correspond to elements of the form $\alpha_h$. By
Lemma \ref{1t1}, for any $h\in\Phi$ there exists a unique edge
$\te$ of $\wtk$ such that $g_{\te}=h$. Take any $1$-path
$\tq=\tq'e\tq''$ from $0$ to $1$ in $\wtk$ containing $\te$. By
the property T2, $\tx$ and $\tq$ are homotopic so there exists a
$(\tq,\tx)$-diagram $\Delta$ over $\wtk$. Now we consider the
product of three factors. The first of them is $\Delta^{-1}$, the
third one is $\Delta$, and the second factor is the atomic diagram
$\ve(\tq')+f+\ve(\tq'')$, where $f$ is the positive leaf of the
form $\te=\te$. By Theorem \ref{aaapc}, the diagram $A_h$
representing $a(\te)$ is obtained from this diagram by changing
all edge labels to $x$. As a result, we get the diagram (an
element of $\gr_1$).

Now let $g\in F$ represent an $(x,x)$-diagram $\Gamma$ over
$\hh_0$. So $g^{-1}\alpha_hg\in\gr_1$ is represented by the diagram
$\Gamma^{-1}A_h\Gamma$. We claim that it is equal to $A_{hg}$ in
$\gr_1$.

The diagram $\Gamma^{-1}A_h\Gamma$ has a unique $(x,x)$-cell so it
can be written as the product of three factors, where the second
factor is an atomic diagram with a positive $(x,x)$-cell. Each
diagram of this form represents a generator of $\aaa$. To find
which one, we must take the first factor and lift it into $\wtk$
according to Lemma \cite[Lemma 8.3]{GuSa02a}. The top edge of the
$(x,x)$-cell will be mapped onto some edge $\bar e$ of $\wtk$.
This will tell us what will be the generator of $\aaa$ we want to
find. In our case lift of the first factor is a product of
$\Gamma^{-1}$ and $\Delta^{-1}$. So the transition function
$g_{\bar e}$ of $\bar e$ is the product of the transition function
that corresponds to $\Delta$, which is $T_{\te}=h$, and the
transition function of $\Gamma$, which is $g$. So the composition
of the corresponding functions is $hg$, as desired. Thus
$g^{-1}\alpha_hg=\alpha_{hg}$.
\endproof

Now we finally have the following structural description of $\gr_1$.

\begin{thm}
\label{strg1} The universal diagram group $\gr_1$ is isomorphic to
the semi-direct product of the group $\aaa$ and the R.\,Thompson's
group $F$. Here $\aaa$ is a partially commutative group generated
by symbols of the form $\alpha_h$, where $h$ runs over $\Phi$, the
set of all continuous piecewise linear functions from $[0,1]$ into
$[0,1]$ with finitely many singularity points, where $(0)f$, $(1)f$
and all the singularity points of $h$ are dyadic rational and all
slopes are of the form $2^k$, $k\in\mathbb{Z}$. Two symbols
$\alpha_{h_1}$ and $\alpha_{h_2}$ commute if and only if the
interiors of the images of $h_1$ and $h_2$ are disjoint. The group
$F$, as the group of piecewise linear functions, is a subset in
$\Phi$. Its action on $\aaa$ is defined as follows:
$g^{-1}\alpha_hg=\alpha_{hg}$.
\end{thm}

\section{Orderability}
\label{s-ord}

Now we are finally in a position to prove the main result of the
paper:

\begin{thm}
\label{ord}
Any diagram group of any directed $2$-complex is totally orderable.
\end{thm}

\proof It is well known that the property of a group to be
orderable is local \cite{KK}. This means that the group is
orderable if and only if all its finitely generated subgroups are
orderable. If we have finitely many diagrams as elements of some
diagram group, then we can find a finite subcomplex in the
original directed $2$-complex such that all these diagrams will be
diagrams over a finite directed $2$-complex, that is, they will be
elements of some countable diagram group $H$. Since $\gr_1$ is
universal \cite[Theorem 5.6]{GuSa02a}, the group $H$ embeds into
$\gr_1$. If we can prove that $\gr_1$ is orderable, then $H$ will
be also orderable. So any diagram group will be orderable.

By Theorem \ref{strg1}, $\gr_1$ is a semi-direct product of $\aaa$ and $F$.
The group $\aaa$ is partially commutative and thus every total order on the
set of generators of $\aaa$ extends to a total order of $\aaa$ by \cite{DuKr92}.
The group $F$ is also orderable \cite{BrSq85}. It is easy to see that to show
that $\gr_1$ is orderable, it suffices to show that the action of $F$ on $\aaa$
respects a total order on $\aaa$.

There exists a well known an natural total order on the set $\Phi$
\cite{BrSq85}. Let $h_1,h_2\in\Phi$. Let $c$ be the least upper
bound of the set $\{\,t\mid h_1\mbox{\ equals\ }h_2\mbox{\ on\
}[0,t]\,\}$. If $h_1\ne h_2$, then $c<1$. Both functions have the
same value at $c$ and they are affine on a small right
neighbourhood of $c$. So they have different right derivatives at
$c$. Denoting these derivatives by $k_1$ and $k_2$, respectively,
we set $h_1<h_2$ whenever $k_1<k_2$ and $h_2<h_1$ whenever
$k_2<k_1$. Since $F$ is a subset in $\Phi$, this induces also the
order on $F$. This order on $F$ is stable under multiplication
from both sides \cite{BrSq85}. We also need to notice that
$h_1<h_2$ implies $h_1g<h_2g$ whenever $h_1,h_2\in\Phi$, $g\in F$.
This is obvious because all functions from $\Phi$ are increasing.

We can introduce a total order on the set of generators of $\aaa$
in the following natural way: $\alpha_{h_1}<\alpha_{h_2}$ if and
only if $h_1<h_2$ in $\Phi$. This order induces a total order on
$\aaa$ as described in \cite{DuKr92}. Let \be{lcs}
\aaa=\aaa_1\supseteq\aaa_2\supseteq\cdots\supseteq\aaa_n\supseteq\cdots
\ee be the lower central series of $\aaa$:
$\aaa_{n+1}=[\aaa_n,\aaa]$ for all $n\ge1$. It is known
\cite{DuKr92} that the intersection of all the $\aaa_n$'s is
trivial for any partially commutative group. Following
\cite{DuKr92}, we define basic commutators of weight $n\ge1$ by
induction and also introduce a total ordering on the set of basic
commutators. By definition, basic commutators of weight $1$ are
the elements $\alpha_h$ ($h\in\Phi$). We already have a total
order on that set. Now let $n>1$ and suppose that we have a
definition of basic commutators of weight $<n$ and also a total
order $<$ on them. Basic commutators of weight $n$ will be bigger
than basic commutators of weight $<n$ with respect to the order we
define. Any basic commutator of weight $n>1$ is a commutator of
the form $[c',c'']$, where

\begin{itemize}
\item $c'$, $c''$ are basic commutators of weight $<n$ and the sum
of their weights is $n-1$,
\item $[c',c'']$ is not equal to 1 in $\aaa$.
\item $c'>c''$ in the order we have,
\item if $c'=[c_1,c_2]$ for some basic commutators $c_1$, $c_2$, then
$c''\ge c_2$.
\end{itemize}

Basic commutators of weight $n$ are ordered lexicographically:
$[c_1,c_2]<[c_3,c_4]$ if and only if either $c_1<c_2$, or
$c_1=c_2$ and $c_3<c_4$.

By \cite{DuKr92}, basic commutators of weight $n$ freely generate
the Abelian group $\aaa_n/\aaa_{n+1}$. These groups can be totally
ordered in a standard way. Say, if $Y$ is the basis of a free
Abelian group totally ordered with $<$, then any nontrivial
element $g$ of the group generated by $Y$ is a (uniquely defined)
product of powers of elements of $Y$ with non-zero exponents:
$g=y_1^{k_1}\cdots y_m^{k_m}$, where $y_1<\cdots<y_m$. By
definition, $g$ is positive whenever the exponent on $y_m$ is
positive. By \cite[Theorem 3.1]{DuKr92}, there exists a unique
total order on $\aaa$ such that the canonical projections
$\pi_n\colon\aaa_n\to\aaa_n/\aaa_{n+1}$ are increasing. To find
out whether an element $g\in\aaa$, $g\ne1$, is positive, it
suffices to find the lowest $n$ such that $g\notin\aaa_{n+1}$.
Since $g\in\aaa_n$, one can project $g$ onto a nontrivial element
$\pi_n(g)\in\aaa_n/\aaa_{n+1}$. One has $g>1$ if and only if
$\pi_n(g)>1$ in the corresponding group.

Notice that the word problem is decidable in all partially
commutative groups (see, for example, \cite{DJ}) so one can
effectively decide if an element of $\aaa$ is positive.

It remains to prove that this order $<$ on $\aaa$ is stable under
the action of $F$. We know that any element $g\in F$ induces a
permutation on the set of generators $\{\,\alpha_h\mid h\in\Phi\,\}$
and preserves the order on this set. Since $[a,b]^g=[a^g,b^g]$ for
any commutator, the action by $g$ also preserves the set of basic
commutators and the order on that set. This implies that the order
on $\aaa$ is preserved as well.
\endproof

We would like to add a few more applications of the results of the
preceding sections. First we have the following surprising result.

\begin{thm}\label{rescount} Every diagram group is residually
countable.
\end{thm}

\proof The proof of \cite[Lemma 5.3]{GuSa02a} shows that any
diagram group embeds into the diagram group $\dd(\hh_\alpha,x)$,
where $\hh_\alpha$ is the expansion of the Dunce hat with
$\nu(x)=\alpha$ (here $\alpha$ is a cardinal number). Any diagram
over $\hh_\alpha$ involves only finitely many leaf cells. So for
any reduced nontrivial diagram $\Delta$ over $\hh_\alpha$, we have
a homomorphism from $\dd(\hh_\alpha,x)$ to the countable group
$\gr_n=\dd(\hh_n,x)$ (here $n$ is the number of leaf cells in
$\Delta$) that collapses all leaves not involved in $\Delta$ to
edges. This homomorphism takes $\Delta$ to itself, so it separates
$\Delta$ from 1.
\endproof

We know that there are simple non-trivial diagram groups: the
derived subgroup of the group $F$ is an example  (see
\cite[Theorem 26]{GuSa99}). Other examples include the derived
subgroups of all generalized Thompson groups $F_n$ (see
\cite{Bro87} or \cite{GuSa02a} for definitions). It is known
\cite{BrGu98} that all generalized Thompson groups are embeddable
into each other. In particular, all of them are subgroups in $F$.
Here is our result.

\begin{thm}
\label{smpl} Any simple subgroup of a diagram group embeds into
R.\,Thomp\-son's group $F$.
\end{thm}

\proof Let $G$ be a simple subgroup of a diagram group. By Theorem
\ref{rescount}, $G$ is countable. Therefore it is embeddable into
$\gr_1$ by \cite[Theorem 5.6]{GuSa02a}. Since $\gr_1$ is an
extension of $\aaa$ by $F$, the group $G$ either embeds into
$\aaa$ or it embeds into $F$. The first option is impossible
because every partially commutative group is residually nilpotent
\cite{DuKr92}. Hence $G$ is a subgroup of $F$.
\endproof

It is natural to ask the following question.

\begin{prob}
\label{pr1} Is it true that any simple subgroup of $F$ is
isomorphic to the derived subgroup of $F_n$ for some $n$?
\end{prob}

One more corollary deals with groups representable by diagrams
that do not contain non-Abelian free subgroups. We know that $F$
is one of such groups \cite{CFP}. The following question is still
open.

\begin{prob}
\label{pr2} Is it true that every group representable by diagrams
and containing no non-Abelian free subgroups embeds into $F$?
\end{prob}

However, the following result is a step in this direction.

\begin{thm}
\label{embf} Let $G$ be a group representable by diagrams that
does not contain non-Abelian free subgroups. Then $G$ is an
extension of some free Abelian group $A$ by a subgroup $B$ of the
group $F$.
\end{thm}

This theorem will follow from the next lemma.

\begin{lm} \label{DJ} Every subgroup of a partially commutative group
is either Abelian or contains a free non-Abelian subgroup.
\end{lm}

\proof Indeed, it suffices to prove the statement for subgroups of
finitely generated partially commutative groups. By \cite{DJ},
every finitely generated partially commutative group has a finite
index subgroup embeddable into a Coxeter group.

By \cite{NV}, every finitely generated subgroup of a Coxeter group
is either virtually Abelian or contains a finite index subgroup
with a non-Abelian free factor. Having a free non-Abelian factor
implies having a free non-Abelian subgroup.

Since partially commutative groups are orderable, they satisfy the
{\em centralizer property} \cite{KK}, that is $[a^n,b]=1$ implies
$[a,b]=1$ for every pair of elements $a,b$. This immediately
implies that every virtually Abelian subgroup of a partially
commutative group is Abelian.
\endproof

Now it is easy to complete the proof of Theorem \ref{embf}.

\proof Let $H$ be a subgroup of a diagram group. Suppose that all
free subgroups of $G$ are Abelian. As above, we can regard $G$ as
a subgroup of $\gr_\alpha=\dd(\hh_\alpha,x)$ for some cardinal
$\alpha$. By Theorem \ref{aaapc}, $\gr_\alpha$ is an extension of
a partially commutative group $\aaa_\alpha$ by $F$. Then $G$ is an
extension of $G\cap \aaa_\alpha$ by a subgroup of $F$. By Lemma
\ref{DJ}, $G\cap \aaa_\alpha$ is Abelian. It remains to apply
\cite[Theorem 16]{GuSa99} that says that every Abelian subgroup of
a diagram group is free Abelian.
\endproof

\vspace{10ex}

\begin{minipage}[t]{2.9in}
\noindent Victor Guba\\
Department of Mathematics\\
Vologda State University\\
guba@uni-vologda.ac.ru
\end{minipage}
\begin{minipage}[t]{2.9in}
\noindent Mark V. Sapir\\
Department of Mathematics\\
Vanderbilt University\\
msapir@math.vanderbilt.edu
\end{minipage}

\end{document}